\documentclass{article}
\usepackage{amssymb,amsmath,latexsym, verbatim}
\usepackage{hyperref}
\usepackage{cite}

\newtheorem{thm}{Theorem}[section]
\newtheorem{lemma}[thm]{Lemma}
\newtheorem{cond}[thm]{Condition}
\newtheorem{defn}[thm]{Definition}
\newtheorem{prop}[thm]{Proposition}
\newtheorem{corollary}[thm]{Corollary}
\newtheorem{remark}[thm]{Remark}
\newtheorem{example}[thm]{Example}

\numberwithin{equation}{section}

\def\qed{{\hfill $\Box$ \bigskip}}

\def\R{{\mathbb R}}

\def\E{{\mathbb E}}

\def\P{{\mathbb P}}
\def\N{{\mathbb N}}
\def\eps{\varepsilon}

\def\pf{\noindent{\bf Proof.} }
\def\beq{\begin{equation}}
\def\eeq{\end{equation}}
\def\lam{\lambda}

\allowdisplaybreaks

\begin{document}

\begin{doublespace}
\title{\Large \bf  Fatou's theorem for subordinate Brownian motions with Gaussian components on $C^{1,1}$ open sets}
\author{Hyunchul Park}

\date{ \today}
\maketitle

\begin{abstract}
We prove Fatou's theorem for nonnegative harmonic functions with respect to killed subordinate Brownian motions with Gaussian components on bounded $C^{1,1}$ open sets $D$. 
We prove that nonnegative harmonic functions with respect to such processes on $D$ converge nontangentially almost everywhere with respect to the surface measure as well as the harmonic measure restricted to the boundary of the domain. 
In order to prove this, we first prove that the harmonic measure restricted to $\partial D$ is mutually absolutely continuous with respect to the surface measure. We also show that tangential convergence fails on the unit ball.
\end{abstract}


\noindent{\bf Key works and phrases}: Fatou theorem, subordinate Brownian motions with Gaussian component, Martin representation, harmonic function

\section{Introduction}
In \cite{Fatou} Fatou showed that bounded classical harmonic functions in the unit disc converge nontangentially almost everywhere. 
The nontangential convergence of harmonic functions is generally called Fatou's theorem. 
Later the Fatou's theorem for diffusion processes is extended into many directions. 
The Fatou theorem is established on more general domains and in \cite{Ai} the author proved the Fatou's theorem for classical harmonic functions on uniform domains. 
The other direction is to establish the Fatou's theorem for more general operators than Laplacian. 
In \cite{CDM} the authors established Fatou's theorem for a family of elliptic operators in the unit ball in $\mathbb{C}^{d}$, $d\geq 2$.
These results deal with the nontangential convergence of  harmonic functions with respect to operators that correspond to diffusion processes. 

The Fatou's theorem is also established for pure jump processes. In \cite{BY1, BY2} the authors showed that regular harmonic functions with respect to stable processes (see \hyperref[defn:harmonic]{\ref{defn:harmonic}} for definition of the regular harmonic functions) converge nontangentially almost everywhere on half-spaces and Lipschitz domains, respectively. 
The Fatou's theorem is established for other jump processes and in \cite{Kim2} the author proved the Fatou's theorem for harmonic functions with respect to censored stable processes on bounded $C^{1,1}$ open sets. 

Recently, there have been many interests about Markov processes that have both diffusion and jump components.
A typical prototype of these processes would be an independent sum of a Brownian motion with a symmetric stable process and their potential theoretical properties have been investigated in \cite{CKS1, CKS2, CKSV}.
In \cite{KSV1} the authors studied subordinate Brownian motions (SBMs) with Gaussian components and established the boundary Harnack principle for harmonic functions with respect to such processes, established sharp two-sided Green function estimates, and identified the Martin boundary with the Euclidean boundary of $C^{1,1}$ open sets. 

In this paper we consider subordinate Brownian motions $X$ with a diffusion component and a quite general jump component (see Section \hyperref[Preliminaries]{\ref{Preliminaries}} for precise definition).
The main goal of this paper is to prove that the Fatou theorem holds true for nonnegative harmonic functions with respect to killed processes $X^{D}$ (Corollary \hyperref[cor:Fatou killed]{\ref{cor:Fatou killed}}) on bounded $C^{1,1}$ open sets $D$. 
We prove that the nontangential convergence occurs almost every point with respect to the surface measure of $\partial D$. Note that this is a very different situation from a case when underlying processes are pure jump processes. For regular harmonic functions with respect to symmetric stable processes the authors in \cite{BY1, BY2} showed that the Fatou theorem for stable processes requires more restrictive conditions than for Brownian motions. For harmonic functions with respect to killed symmetric stable processes, it is  proved that a certain harmonic function is comparable to $\delta_{D}(x)^{\frac{\alpha}{2}-1}$, which can not converge as the point $x$ approaches a point in $\partial D$ in \cite[Equation 11]{BD}. 

We prove the Fatou theorem using both analytical techniques and probabilistic techniques that are similar to \cite{Kim2} or \cite{Kim}. 
However, this process is not straightforward. In the probabilistic techniques to establish the Fatou theorem for censored stable processes in \cite{Kim2} or the relative Fatou theorem for symmetric stable processes in \cite{Kim}, the author identified the probabilistic martingale convergence of nonnegative harmonic functions with analytical nontangential convergence and an oscillation estimate for harmonic functions on balls of different radii played an important role (\cite[Proposition 3.9]{Kim2} and \cite[Proposition 3.11]{Kim}). The oscillation estimate for harmonic functions on balls with respect to those processes is quite straightforward due to explicit expressions for Poisson kernels for balls for stable processes. Such explicit expressions for the Poisson kernels are not available anymore for general subordinate Brownian motions with Gaussian components. To overcome this difficulty, we first establish relative Fatou theorem for harmonic functions with respect to killed processes (Theorem \ref{thm:RFT}) and we use the relative Fatou theorem to identify the probabilistic convergence with the analytic nontangential convergence (Proposition \hyperref[p:2.12]{\ref{p:2.12}}). 

This paper is organized as follows. In Section \hyperref[Preliminaries]{\ref{Preliminaries}} we define the subordinate Brownian motions with Gaussian components and state a few properties about them. In Section \hyperref[Relative Fatou theorem]{\ref{Relative Fatou theorem}} we establish the relative Fatou's theorem for harmonic functions with respect to $X^{D}$, which asserts the existence of nontangential limit of the ratio of harmonic functions with respect to killed processes. In Section \hyperref[Harmonic measure]{\ref{Harmonic measure}} we prove the Martin measure of harmonic function $F(x):=\P_{x}(X_{\tau_{D}}\in\partial D)$, which represents the probability that the processes exit the domain through the boundary, is absolutely continuous with respect to the surface measure (Theorem \ref{thm:Rep harmonic boundary}) and establish the Fatou theorem for $X^{D}$ (Corollary \ref{cor:Fatou killed}). 
In Section \hyperref[Integral representation of harmonic functions and Fatou theorem]{\ref{Main result}} we establish an integral representation theorem (Theorem \hyperref[thm:IntRep]{\ref{thm:IntRep}}) for harmonic functions with respect to $X$. 
We also show that our result is best possible by showing that tangential convergence of harmonic functions on the unit ball can fail. 

In this paper the upper case constants $\Lambda,R_{0},R_{1}, R_{2}, C_{1},C_{2},C_{3},C_{4}$ will be fixed. The value of lower case constants $\eps$, $\delta$, $\eta$, $c$ or $c_{1}, c_{2}, \cdots$ will not be important and may change from line to line.

\section{Preliminaries}\label{Preliminaries}
In this section we define subordinate Brownian motions with Gaussian components and state some properties about them.
Recall that a subordinator $S=(S_{t}, t\geq 0)$ is an one-dimensional L\'evy process taking values on $[0,\infty)$ with increasing sample paths. A subordinator $S$ can be characterized by its Laplace exponent $\phi$ through the relation
$$
\E[e^{-\lam S_{t}}]=e^{-t\phi(\lam)}, \quad t>0, \lam>0.
$$
A smooth function $\phi : (0, \infty) \to [0, \infty)$ is called a Bernstein function if $(-1)^n D^n \phi \le 0$ for every positive integer $n$.
The Laplace exponent $\phi$ of a subordinator is a Bernstein function with $\phi(0+)=0$ and can be written as
$$
\phi(\lam)=b\lam+\int_{(0,\infty)}(1-e^{-\lam t})m(dt), \quad \lam >0,
$$
where $b\geq 0$ and $m$ is a measure on $(0,\infty)$ satisfying $\int_{(0,\infty)}(1\wedge t)m(dt)<\infty$. $m$ is called the L\'evy measure of $\phi$. In this paper we will assume that $b>0$ in order to have a nontrivial diffusion part for subordinate Brownian motions.
Without lose of generality we assume $b=1$.

Suppose that $W=(W_{t}: t\geq 0)$ is a $d$-dimensional Brownian motion and $S=(S_{t}: t\geq 0)$ is a subordinator
with Laplace exponent $\phi$, which is independent of $W$.
The process $X=(X_{t}:t\geq 0)$ defined by $X_{t}=W(S_{t})$ is called a subordinate Brownian motion and its infinitesimal generator is given by $\phi(\Delta):=-\phi(-\Delta)$, which can be constructed via Bochner's functional calculus. On $C_{b}^{2}(\R^{d})$ (the collection of $C^{2}$
functions in $\R^{d}$ which, along with partial derivatives up to order 2, are bounded),
$\phi(\Delta)$ is an integro-differential operator of the type
$$
\Delta f(x)+\int_{\R^{d}}\left(f(x+y)-f(x)-\nabla f(x)\cdot y 1_{\{|y|\leq 1\}}\right)J(dy),
$$
where the measure $J$ has the form $J(dy)=j(|y|)dy$ with $j:(0,\infty)\rightarrow (0,\infty)$ given by
$$
j(r)=\int_{0}^{\infty}(4\pi t)^{-d/2}e^{-r^{2}/(4t)}m(dt).
$$

Throughout this paper we will impose two conditions on $\phi$ and $m$.
\begin{cond}\label{con:con}
\begin{enumerate}
\item
The Laplace exponent $\phi$ of $S$ is a completely Bernstein function. That is, the L\'evy measure $m$ has a completely monotone density
(i.e., $m (dt) = m(t)dt$ and $(-1)^n D^n m \ge 0$ for every non-negative integer $n$).
\item For any $K>0$, there exists $c=c(K)>1$ such that
$$
m(r)\leq cm(2r)
\quad \text{for }\,\,  r\in (0,K).
$$
\end{enumerate}
\end{cond}
Note that Condition \hyperref[con:con]{\ref{con:con}} is the main assumption imposed in \cite{KSV1}.

There are many important subordinators that satisfy Condition \hyperref[con:con]{\ref{con:con}} and we list some of most important examples. 
\begin{example}
\begin{enumerate}
\item A function $\ell(x)$ is slowly varying at $\infty$ if $\displaystyle\lim_{x\rightarrow\infty}\frac{\ell(\lam x)}{\ell(x)}=1$ for all $\lam>0$. Let $\phi(\lam)$ be a complete Bernstein function which satisfies
$$
\lam +c_{1}\lam^{\alpha/2}\ell(\lam)\leq \phi(\lam)\leq \lam +c_{2}\lam^{\alpha/2}\ell(\lam),
$$
for some constants $0<c_{1}\leq c_{2}<\infty$, $0<\alpha<2$, and $\ell(\lam)$ is slowly varying at $\infty$. It follows from \cite[Theorem 2.10]{KSV2} that Conditions \hyperref[con:con]{\ref{con:con}} are satisfied for those processes. In particular these classes contain the sum of Brownian motions and symmetric stable processes, relativistic stable processes with mass $m$, and mixed stable processes and the corresponding $\phi(\lam)$ are given by $\phi(\lam)=\lam +\lam^{\alpha/2}$, $\phi(\lam)=\lam+ \left((m^{2/\alpha}+\lam)^{\alpha/2}-m\right)$, $\phi(\lam)=\lam+\lam^{\alpha/2}+\lam^{\beta/2}$, $0<\beta<\alpha<2$, respectively. 

\item Geometric stable subordinator 

Let $\phi(\lam)=\lam + \ln(1+\lam^{\alpha/2})$, $0<\alpha\leq 2$.  From \cite[Theorem 2.4]{SSV} Conditions \hyperref[con:con]{\ref{con:con}} are satisfied. 
Note that when $\alpha=2$ it corresponds to the sum of Brownian motions and Gamma processes and the corresponding L\'evy density is given by $m(t)=\frac{e^{-t}}{t}$.
\end{enumerate}
\end{example}

For any open set $D\subset \R^{d}$, $\tau_{D}:=\inf\{t>0 : X_{t}\notin D\}$ denotes the first exit time from $D$ by $X$.
We will use $X^{D}$ to denote the process defined by $X_{t}^{D}(\omega)=X_{t}(\omega)$
if $t<\tau_{D}(\omega)$ and
$X_{t}^{D}(\omega)=\partial$ if $t\geq \tau_{D}(\omega)$, where $\partial$ is a cemetery point.
It is well known that $X^{D}$ is a strong Markov process with state space $D\cup \{\partial\}$.
For any function $u(x)$ defined on $D$ we extend it to $D\cup \{\partial\}$ by letting $u(\partial)=0$.
It follows from \cite[Chapter 6]{BBKRSV} that the process $X$ has a transition density $p(t,x,y)$ which is jointly continuous. Using this and the strong Markov property, one can easily check that
$$
p_{D}(t,x,y):=p(t,x,y)-\E_{x}[p(t-\tau_{D},X_{\tau_{D}},y); t>\tau_{D}], \quad x,y\in D
$$
is continuous and is a transition density of $X^{D}$.
For any bounded open set $D\subset \R^{d}$, we will use $G_{D}(x,y)$ to denote the Green function of $X^{D}$, i.e.,
$$
G_{D}(x,y):=\int_{0}^{\infty}p_{D}(t,x,y)dt, \quad x,y\in D.
$$
Note that $G_{D}(x,y)$ is continuous on $\{(x,y)\in D\times D : x\neq y\}$.

The L\'evy density is given by $J(x,y)=j(|x-y|)$, $x,y\in \R^{d}$ and it determines the L\'evy system for $X$, which describes the jumps of the process $X$:
For any nonnegative measurable function $F$ on $\R_{+}\times\R^{d}\times \R^{d}$ with $F(s,x,x)=0$ for all $s>0$ and $x\in\R^{d}$, and stopping time $T$
with respect to $\{\mathcal{F}_t: t\ge 0\}$,
$$
\E_{x}\left[\sum_{s\leq T}F(s,X_{s-},X_{s})\right]=\E_{x}\left[\int_{0}^{T}\left(\int_{\R^{d}}F(s,X_{s},y)J(X_{s},y)dy\right)ds\right].
$$
Using the L\'evy system, we know that for any nonnegative function $f\geq 0$ and every bounded open set $D$ we have
\beq\label{eqn:Levy system}
\E_{x}\left[f(X_{\tau_{D}}), X_{\tau_{D^{-}}}\neq X_{\tau_{D}}\right]=\int_{\overline{D}^{c}}\int_{D}G_{D}(x,y)J(y,z)dy f(z)dz,\quad x\in D.
\eeq
We define $K_{D}(x,z)=\int_{D}G_{D}(x,y)J(y,z)dy$ and \hyperref[eqn:Levy system]{\eqref{eqn:Levy system}} can be written as
\beq\label{eqn:Levy system2}
\E_{x}\left[f(X_{\tau_{D}}), X_{\tau_{D^{-}}}\neq X_{\tau_{D}}\right]=\int_{\overline{D}^{c}}K_{D}(x,z)f(z)dz,\quad x\in D.
\eeq

Now we state the definition of harmonic functions.
\begin{defn}\label{defn:harmonic}
\begin{enumerate}
\item A function $u:\R^{d}\rightarrow [0,\infty)$ is said to be
harmonic in an open set $D\subset\R^{d}$ with respect to $X$ if for every open set $B$ whose closure is a compact subset of $D$,
$$
u(x)=\E_{x}\left[u(X_{\tau_{B}})\right]\quad \text{for every } x\in B.
$$
\item A function $u:\R^{d}\rightarrow [0,\infty)$ is said to be
regular harmonic in $D$ with respect to $X$ if 
$$
u(x)=\E_{x}\left[u(X_{\tau_{D}})\right]\quad \text{for every } x\in D.
$$
\item A function $u:\R^{d}\rightarrow [0,\infty)$ is said to be
harmonic with respect to $X^{D}$ if for every open set $B$ whose closure is a compact subset of $D$,
$$
u(x)=\E_{x}\left[u(X^{D}_{\tau_{B}})\right]= \E_{x}\left[u(X_{\tau_{B}}), \tau_{B}<\tau_{D}\right]\quad \text{for every } x\in B.
$$
\end{enumerate}
\end{defn}
Note that it follows from strong Markov property that every regular harmonic function is harmonic.

The following Harnack principle is proved in \cite[Proposition 2.2]{KSV1}.
\begin{prop}\label{prop:HP}
There exists a constant $c>0$ such that for any $r\in(0,1]$, $x_{0}\in \R^{d}$ and any function $f$
which is nonnegative in $\R^{d}$ and harmonic in $B(x_{0},r)$ with respect to $X$ we have
$$
f(x)\leq cf(y) \quad \text{for all } x,y\in B(x_{0},r/2).
$$
\end{prop}

Recall that an open set $D$ in $\R^{d}$ is said to be a (uniform) $C^{1,1}$ open set
if there are (localization radius) $R_{0}>0$
and $\Lambda_{0}$ such that for every $z\in\partial D$ there exist
a $C^{1,1}$ function $\psi=\psi_{z}:\R^{d}\rightarrow \R$ satisfying $\psi(0, \cdots, 0)=0$,
$\nabla\psi(0)=(0,\cdots,0)$, $|\nabla\psi(x)-\nabla\psi(y)|\leq \Lambda_{0}|x-y|$, and an orthonormal coordinate system
$CS_{z}:y=(y_{1},\cdots,y_{d-1},y_{d}):=(\tilde{y},y_{d})$ with its origin at $z$
such that $B(z,R_{0})\cap D =\{y= (\tilde{y},y_{d})\in B(0,R_{0}) \text{ in } CS_{z} : y_{d}>\psi(\tilde{y})\}$.
In this paper we will call the pair $(R_{0},\Lambda_{0})$ the characteristics of the $C^{1,1}$ open set $D$.

We state the result about the Martin boundary of a bounded $C^{1,1}$ open set $D$ with respect to $X^{D}$. 
For the definition and its basic properties of the Martin boundary we refer readers to \cite{KW}.
Fix $x_{0}\in D$ and define
$$
M_{D}(x,y):=\frac{G_{D}(x,y)}{G_{D}(x_{0},y)},\quad x,y\in D, \,\, y\neq x, x_{0}.
$$
A positive harmonic function $f$ with respect to $X^{D}$ is called minimal if, whenever $g$ is a positive harmonic function with respect to $X^{D}$ with $g\leq f$,
one must have $f=cg$ for some positive constant $c$.
Now we recall the identification of the Martin boundary of bounded $C^{1,1}$ open sets $D$ with respect to killed processes $X^{D}$ with the Euclidean
boundary in \cite{KSV1}.
\begin{thm}\emph{(\cite[Theorem 1.5]{KSV1})}\label{thm:Martin}
Suppose that $D$ is a bounded $C^{1,1}$ open set in $\R^{d}$. For every $z\in \partial D$, there exists $M_{D}(x,z):=\displaystyle\lim_{y\rightarrow z}M_{D}(x,y)$. Furthermore, for every $z\in \partial D$, $M_{D}(\cdot,z)$ is a minimal harmonic function
with respect to $X^{D}$ and $M_{D}(\cdot,z_{1})\neq M_{D}(\cdot,z_{2})$ for $z_{1},z_{2}\in\partial D$, $z_1 \neq z_2$.
Thus the minimal Martin boundary of $D$ can be identified with the Euclidean boundary.
\end{thm}
Thus by the general theory of Martin boundary representation in \cite{KW} and Theorem \hyperref[thm:Martin]{\ref{thm:Martin}}, we conclude that for every harmonic function
$u\geq 0$ with respect to $X^{D}$, there exists a unique finite measure $\mu$ supported on $\partial D$ such that $u(x)=\int_{\partial D}M_{D}(x,z)\mu(dz)$.
$\mu$ is called the Martin measure of $u$.
\newline

Finally we observe that
the Martin kernel $M_{D}(x,z)$ has the following two-sided estimates.
Let $\delta_{D}(x)=\inf\{|x-z|: z\in D^{c}\}$ be the distance of $x$ from $D^{c}$.
\begin{prop}\label{prop:MK}
Suppose that $D$ is a bounded $C^{1,1}$ open set in $\R^{d}$, $d\geq 2$. Then there exist constants $c_{1}=c_{1}(d,D,\phi)$ and $c_{2}=c_{2}(d,D,\phi)$ such that
\beq\label{eqn:Martin}
c_{1}\frac{\delta_{D}(x)}{|x-z|^{d}} \leq M_{D}(x,z) \leq c_{2}\frac{\delta_{D}(x)}{|x-z|^{d}}, \quad x\in D, \,\,z\in \partial D.
\eeq
\end{prop}
\pf
Let $$
g_{D}(x,y) :=\begin{cases}
\frac{1}{|x-y|^{d-2}}\left(1\wedge \frac{\delta_{D}(x)\delta_{D}(y)}{|x-y|^{2}}\right)&\mbox{when  } d\geq 3,\\
\log \left(1+\frac{\delta_{D}(x)\delta_{D}(y)}{|x-y|^{2}}\right) &\mbox{when } d=2.
\end{cases}
$$
Then it follows from \cite[Theorem 1.4]{KSV1} there exists $c_{1}=c_{1}(d,D,\phi)$ and $c_{2}=c_{2}(d,D,\phi)$ such that
\beq\label{eqn:Green}
c_{1}g_{D}(x,y)\leq G_{D}(x,y) \leq c_{2}g_{D}(x,y).
\eeq
From Theorem \hyperref[thm:Martin]{\ref{thm:Martin}} the martin kernel $M_{D}(x,z)$ can be obtained by $M_{D}(x,z)=\displaystyle\lim_{y\rightarrow z}\frac{G_{D}(x,y)}{G_{D}(x_{0},z)}$. Now from \eqref{eqn:Green} we immediately get the assertion of the proposition.
\qed

\section{Relative Fatou theorem for harmonic functions with respect to $X^{D}$}\label{Relative Fatou theorem}
Throughout this section we assume that $D$ is a bounded $C^{1,1}$ open set in $\R^{d}$, $d\geq 2$ with the characteristics $(R_{0},\Lambda_{0})$.
In this section we prove relative Fatou's theorem for nonnegative harmonic functions $u$ and $v$ with respect to $X^{D}$.
For any finite and nonnegative measure $\mu$ supported on $\partial D$ we define
$$
M_{D}\mu(x):=\int_{\partial D}M_{D}(x,z)\mu(dz), \quad x\in D.
$$
Since $M_{D}(\cdot,z)$ is harmonic with respect to $X^{D}$ for $z\in \partial D$
(see Theorem \hyperref[thm:Martin]{\ref{thm:Martin}}),
it is easy to see
that $M_{D}\mu(x)$ is nonnegative and harmonic with respect to $X^{D}$.

Now we define Stolz open sets. For $z\in\partial D$ and $\beta>1$, let
$$
A_{z}^{\beta}=\{x\in D : \delta_{D}(x)<R_{0} \text{ and } |x-z|<\beta\delta_{D}(x)\}.
$$
We say $x$ approaches $z$ nontangentially if $x\rightarrow z$ and $x\in A_{z}^{\beta}$ for some $\beta>1$.

It is well known that $C^{1,1}$ open sets satisfy uniform interior and exterior ball property with some radius of $R$ (see \cite[Lemma 2.2]{AKSZ}). By decreasing $R_{0}$ in the definition of $C^{1,1}$ open sets if necessary, we may assume $R=R_{0}$. In particular $C^{1,1}$ open sets are $\kappa$-fat open sets with $\kappa=\frac{R_{0}}{2}$ (see \cite{Kim} for the definition of $\kappa$-fat open set). It follows from \cite[Lemma 3.9]{Kim} that for any $z\in\partial D$ and $\beta>\frac{1-\kappa}{\kappa}$
$A_{z}^{\beta}\neq \emptyset$ and there exists a sequence $\{y_{k}\}\subset A_{z}^{\beta}$ such that $\lim_{k\rightarrow\infty} y_{k}=z$.
From now on, we will always assume this condition for $\beta$ so that $A_{z}^{\beta}\neq \emptyset$ for all $z\in\partial D$.

Recall the following property of the surface measure $\sigma$, called \textit{Ahlfors regular condition} (see \cite[page 992]{MR}):
there exist constants $R_{1}=R_{1}(D,d)$, $C_{1}=C_{1}(D,d)$ and $C_{2}=C_{2}(D,d)$ such that for every $z\in\partial D$ and $r\leq R_{1}$
\begin{eqnarray}\label{eqn:Ahlfors}
C_{1}r^{d-1}&\leq& C_{1} \sigma(\partial D\cap (B(z,r)\setminus B(z,r/2)))\leq \sigma(\partial D \cap B(z,r))\nonumber\\
&\leq& C_{2} \sigma(\partial D\cap (B(z,r)\setminus B(z,r/2)))\leq C_{2}r^{d-1}.
\end{eqnarray}
The next lemma is similar to \cite[Lemma 4.4]{MR}. Since we are working on $C^{1,1}$ open sets, the proof is simpler.

\begin{lemma}\label{lemma:positive}
Let $v(x)=M_{D}\nu(x)$, where $\nu$ is a finite and nonnegative measure on $\partial D$.
For $\nu$-almost every point $z\in \partial D$, we have
$$
\liminf_{A_{z}^{\beta}\ni x\rightarrow z\in\partial D}v(x)>0.
$$
In particular $\nu$-almost every point $z\in \partial D$, we have
\beq\label{eqn:cond1}
\lim_{A_{z}^{\beta}\ni x\rightarrow z}\frac{\delta_{D}(x)}{v(x)}=0.
\eeq
\end{lemma}
\pf
If $x\rightarrow z\in\partial D$ nontangentially, there exists a constant $\beta>0$ such that
$$
\delta_{D}(x)\,\leq \,|x-z| \leq \beta\delta_{D}(x).
$$
Take $x\in D$ such that $|x-z|<R_{1}$ and take $w\in B(z,|x-z|)\cap\partial D$.
Then $|x-w|\leq |x-z|+|z-w|\leq 2|x-z|$ so that we obtain $M_{D}(x,w)\geq c_1 M_{D}(x,z)$ by \hyperref[eqn:Martin]{\eqref{eqn:Martin}}.
This implies
\begin{eqnarray*}
v(x)&\geq& \int_{\partial D \cap B(z,|x-z|)}M_{D}(x,w)\nu(dw) \,\geq\, c_1\, M_{D}(x,z) \, \nu(B(z,|x-z|) \cap \partial D) \\
 &\geq& c_2 \, \frac{\delta_{D}(x)}{|x-z|^{d}} \,\nu(B(z,|x-z|)\cap \partial D) \,\geq\,  c_3\,  \frac{\nu(B(z,|x-z|)\cap \partial D)}{|x-z|^{d-1}}.
\end{eqnarray*}
By \hyperref[eqn:Ahlfors]{\eqref{eqn:Ahlfors}} we have $\sigma(B(z,|x-z|\cap\partial D))\geq c_4|x-z|^{d-1}$ for some constant $c_{4}(D,d)$. Hence we have
$$
\frac{\sigma(B(z,|x-z|)\cap \partial D)}{\nu(B(z,|x-z|)\cap \partial D)}\,\geq\, c_5\, \frac{1}{v(x)},
$$
and by \cite[Theorem 5]{Be} the symmetric derivative
$$
\limsup_{x\rightarrow z}\frac{\sigma(B(z,|x-z|)\cap \partial D)}{\nu(B(z,|x-z|)\cap \partial D)}
$$
is finite $\nu$-almost every point $z\in\partial D$.
\qed

The next lemma is an analogue of \cite[Lemma 4.3]{MR}.
\begin{lemma}\label{lemma:con1}
Let $z\in\partial D$ and $v$ be a nonnegative harmonic function with respect to $X^{D}$ with Martin measure $\nu$.
Suppose that $\mu$ is a nonnegative finite measure on $\partial D$.
If $\displaystyle\lim_{x\rightarrow z}\frac{\delta_{D}(x)}{v(x)}=0$, then for every $\eta>0$ we have
$$
\lim_{x\rightarrow z}\frac{\int_{\partial D \cap \{|z-w|\geq \eta \}}M_{D}(x,w)\mu(dw)}{v(x)}=0.
$$
If we assume
$\displaystyle\lim_{x\rightarrow z}\frac{\delta_{D}(x)}{v(x)}=0$ nontangentially,
then the above limit also need be taken nontangentially.
\end{lemma}
\pf
If $|z-w|\geq\eta$ and $|x-z|\leq \eta/2$, then $|x-w|\geq \eta/2$.
Thus from \hyperref[eqn:Martin]{\eqref{eqn:Martin}} we have
\begin{eqnarray*}
\int_{\partial D \cap \{|z-w|\geq\eta\} }M_{D}(x,w)\mu(dw) &\leq& c\int_{\partial D\cap\{|z-w|\geq \eta\}}\frac{\delta_{D}(x)}{|x-w|^{d}}\mu(dw)\\
&\leq&c \, \eta^{-d}\,\delta_{D}(x)\,\mu(\partial D).
\end{eqnarray*}
Hence we have
$$
\frac{\int_{\partial D \cap \{|z-w|\geq \eta\}}M_{D}(x,w)\mu(dw)}{v(x)}\,\leq \, c \, \frac{\mu(\partial D)}{\eta^{d}}
\frac{\delta_{D}(x)}{v(x)}\rightarrow 0
$$
as $x\rightarrow z$.
\qed

\begin{remark}
Note that the condition $\displaystyle\lim_{x\rightarrow z}\frac{\delta_{D}(x)}{v(x)}=0$
cannot be omitted. To see this,
take any points $z,Q\in \partial D$ with $z\neq Q$. Let
$\mu=\nu=\delta_{\{z\}}$ be Dirac measures at $z\in \partial D$, $v(x)=M_{D}\nu(x)=M_{D}(x,z)$, and $\eta=|z-Q|/2$. Then from \hyperref[eqn:Martin]{\eqref{eqn:Martin}},
$\displaystyle\liminf_{x\rightarrow Q}\frac{\delta_{D}(x)}{v(x)} \geq c|z-Q|^{d}>0$.
Clearly
$\frac{\int_{\partial D \cap \{|Q-w|\geq \eta\}}M_{D}(x,w)\mu(dw)}{v(x)}=1$ for any $x\in D$.
\end{remark}

Suppose that $\mu$ and $\nu$ are two measures supported on $\partial D$.
It follows from the Lebesgue-Radon-Nikodym theorem (\cite[Theorem 3.8]{Fo}) 
there exists $\mu_{s}$ singular to $\nu$ and $f\in L^{1}(\partial D,\nu)$ such that $d\mu=fd\nu+d\mu_{s}$. Such a decomposition is called the Lebesgue decomposition.
Consider all points $z\in\partial D$ for which
\beq\label{eqn:Leb}
\lim_{r\rightarrow 0}\frac{\int_{B(z,r)\cap \partial D}\left(|f(w)-f(z)|\nu(dw)+\mu_{s}(dw)\right)}{\nu(B(z,r)\cap\partial D)}=0.
\eeq
It is well-known that $\nu$-a.e. $z\in \partial D$ \hyperref[eqn:Leb]{\eqref{eqn:Leb}} holds true (for example, see \cite[Theorem 3.20 and 3.22]{Fo}).

The next lemma is the nontangential maximal inequality that is analogous to
\cite[Lemma 4.5]{MR}.
\begin{lemma}\label{lemma:maximal}
Suppose that $\mu$ and $\nu$ are nonnegative finite measure on $\partial D$. For any $x\in D$ and $z\in\partial D$ such that $|x-z|\leq t\delta_{D}(x)$ there exist constants 
$c_{1}=c_{1}(d,D,\phi,t)$ and $c_{2}=c_{2}(d,D,\phi,t)$ such that 
$$
c_{1}\,\inf_{r>0}\frac{\mu(B(z,r)\cap\partial D)}{\nu(B(z,r)\cap\partial D)}\,\leq\,\frac{\int_{\partial D}M_{D}(x,w)\mu(dw)}{\int_{\partial D}M_{D}(x,w)\nu(dw)}\,\leq\, c_{2}\,
\sup_{r>0}\frac{\mu(B(z,r)\cap\partial D)}{\nu(B(z,r)\cap\partial D)}.
$$
\end{lemma}
\pf
The proof is similar to \cite[Lemma 4.5]{MR} but we provide the details for the reader's convenience.
Define $B_{n}=B(z,2^{n}|x-z|)\cap \partial D$ for $n\geq 0$ and $A_{0}:=B_{0}$ and $A_{n}=B_{n}\setminus B_{n-1}$ for $n\geq 1$.
Suppose that $w\in B_{1}$. Then $|x-w|\leq |x-z|+|z-w|\leq 3|x-z|$ and $|x-w|\geq \delta_{D}(x)\geq \frac{|x-z|}{t}$.
Hence it follows from \hyperref[eqn:Martin]{\eqref{eqn:Martin}} there exist $c_{3}(d,D,\phi)$ and $c_{4}(d,D,\phi)$ such that 
$$
M_{D}(x,w)\geq c_{3}\frac{\delta_{D}(x)}{|x-w|^{d}}\geq c_{3}\frac{\delta_{D}(x)}{3^{d}|x-z|^{d}},
$$
and
$$
M_{D}(x,w)\leq c_{4}\frac{\delta_{D}(x)}{|x-w|^{d}}\leq c_{4}\frac{t^{d}\delta_{D}(x)}{|x-z|^{d}}.
$$
Hence for any $w,w' \in B_{1}$ we have
$$
M_{D}(x,w)\leq c_{4}t^{d}\frac{\delta_{D}(x)}{|x-z|^{d}}=\frac{c_{4}3^{d}t^{d}}{c_{3}}\frac{c_{3}\delta_{D}(x)}{3^{d}|x-z|^{d}}\leq \frac{c_{4}3^{d}t^{d}}{c_{3}} M_{D}(x,w').
$$
Suppose that $w\in A_{n}$, $n\geq 2$. Then $|x-w|\leq |x-z|+|z-w|\leq (2^{n}+1)|x-z|\leq 2^{n+1}|x-z|$ and 
$|x-w|\geq |w-z|-|x-z|\geq (2^{n-1}-1)|x-z|\geq 2^{n-2}|x-z|$. Hence from \hyperref[eqn:Martin]{\eqref{eqn:Martin}} we have
$$
M_{D}(x,w)\geq c_{3}\frac{\delta_{D}(x)}{|x-w|^{d}}\geq \frac{c_{3}\delta_{D}(x)}{(2^{n+1})^{d}|x-z|^{d}},
$$
and
$$
M_{D}(x,w)\leq c_{4}\frac{\delta_{D}(x)}{|x-w|^{d}}\leq \frac{c_{4}\delta_{D}(x)}{(2^{n-2})^{d}|x-z|^{d}}.
$$
Hence for $w,w'\in A_{n}$, $n\geq 2$ we have
$$
M_{D}(x,w)\leq \frac{c_{4}\delta_{D}(x)}{(2^{n-2})^{d}|x-z|^{d}}=\frac{c_{4}2^{3d}}{c_{3}}\frac{c_{3}\delta_{D}(x)}{(2^{n+1})^{d}|x-z|^{d}}\leq \frac{c_{4}2^{3d}}{c_{3}}M_{D}(x,w').
$$
Set $c_{5}:=\text{max}\left(\frac{c_{4}3^{d}t^{d}}{c_{3}},  \frac{c_{4}2^{3d}}{c_{3}}\right)$. Then we have for any $w,w'\in A_{n}$, $n\geq 0$
\beq\label{eqn:Martin comparable}
M_{D}(x,w)\leq c_{5}M_{D}(x,w').
\eeq

Set $a_{n}:=\sup_{w\in A_{n}}M_{D}(x,w)$ and $b_{n}:=\sup_{k\geq n}a_{k}$ for $n\geq 0$. Clearly $b_{n}\geq a_{n}$ for $n\geq 0$.
Suppose that $w\in A_{n}$ and $w'\in A_{k}$ with $k\geq n+1$. Then $|x-w|\leq 2^{n}|x-z|\leq 2^{k-1}|x-z|\leq |x-w'|$. Hence from \hyperref[eqn:Martin]{\eqref{eqn:Martin}} there exists a constant 
$c_{6}=c_{6}(d,D,\phi)>1$ such that 
$$
M_{D}(x,w')\leq c_{4}\frac{\delta_{D}(x)}{|x-w'|^{d}}\leq c_{4}\frac{\delta_{D}(x)}{|x-w|^{d}}\leq c_{6}M_{D}(x,w).
$$
Hence we have $b_{n}\leq c_{6}a_{n}$ for all $n\geq 0$.

Since $D$ is bounded there exists $k_{0}\in \mathbb{N}$ such that $\partial D\subset \cup_{n=0}^{k_{0}}A_{n}$. Hence it follows from \hyperref[eqn:Martin comparable]{\eqref{eqn:Martin comparable}} and from the fact that $b_{n}\leq c_{6}a_{n}$ we have
\begin{eqnarray*}
&&\int_{\partial D}M_{D}(x,w)\mu(dw)\\
&=&\sum_{n=0}^{k_{0}}\int_{A_{n}}M_{D}(x,w)\mu(dw)\\
&\leq&\sum_{n=0}^{k_{0}}a_{n}\mu(A_{n})\leq\sum_{n=0}^{k_{0}}b_{n}\mu(A_{n})\\
&\leq&b_{0}\mu(B_{0})+\sum_{n=1}^{k_{0}}b_{n}\left(\mu(B_{n})-\mu(B_{n-1})\right)\\
&\leq&\sum_{n=0}^{k_{0}-1}(b_{n}-b_{n+1})\mu(B_{n})+b_{k_{0}}\mu(B_{k_{0}})\\
&\leq&\sum_{n=0}^{k_{0}-1}(b_{n}-b_{n+1})\frac{\mu(B_{n})}{\nu(B_{n})}\nu(B_{n})+b_{k_{0}}\frac{\mu(B_{k_{0}})}{\nu(B_{k_{0}})}\nu(B_{k_{0}})\\
&\leq&\sup_{r>0}\frac{\mu(B(z,r)\cap \partial D)}{\nu(B(z,r)\cap \partial D)}\left(\sum_{n=0}^{k_{0}-1}(b_{n}-b_{n+1})\nu(B_{n})+b_{k_{0}}\nu(B_{k_{0}})\right)\\
&=&\sup_{r>0}\frac{\mu(B(z,r)\cap \partial D)}{\nu(B(z,r)\cap \partial D)}\sum_{n=0}^{k_{0}}b_{n}\nu(A_{n})\\
&\leq&c_{6}\sup_{r>0}\frac{\mu(B(z,r)\cap \partial D)}{\nu(B(z,r)\cap \partial D)}\sum_{n=0}^{k_{0}}a_{n}\nu(A_{n})\\
&\leq&c_{5}c_{6}\sup_{r>0}\frac{\mu(B(z,r)\cap \partial D)}{\nu(B(z,r)\cap \partial D)}\sum_{n=0}^{k_{0}}\int_{A_{n}}M_{D}(x,w)\nu(dw)\\
&=&c_{5}c_{6}\sup_{r>0}\frac{\mu(B(z,r)\cap \partial D)}{\nu(B(z,r)\cap \partial D)}\int_{\partial D}M_{D}(x,w)\nu(dw).
\end{eqnarray*}
Now set $c_{2}:=c_{5}c_{6}$. The opposite inequality can be proved in a similar way and this proves the assertion of the lemma.
\qed

Now we state the main theorem of this section.
\begin{thm}\label{thm:RFT}
Let $u,v$ be nonnegative and harmonic functions with respect to $X^{D}$.
Let $u(x)=\int_{\partial D}M_{D}(x,w)\mu(dw)$ and $v(x)=\int_{\partial D}M_{D}(x,w)\nu(dw)$, where $\mu$ and $\nu$ are nonnegative and finite measures on $\partial D$. 
Let $d\mu=fd\nu+d\mu_{s}$ be Lebesgue decomposition of $\mu$ with respect to
$\nu$.
Then for $\nu$-almost every point $z\in\partial D$ we have
$$
\lim_{x\rightarrow z}\frac{u(x)}{v(x)}=f(z)
$$
as $x\rightarrow z$ nontangentially.
More precisely, the convergence holds for every $z\in \partial D$ satisfying \hyperref[eqn:Leb]{\eqref{eqn:Leb}} and
$\displaystyle\lim_{x\rightarrow z}\frac{\delta_{D}(x)}{v(x)}=0$ as $x\to z$ nontangentially.
\end{thm}
\pf
The proof is similar to
\cite[Theorem 4.2]{MR} but we provide the details for the reader's convenience.
Fix a point $z\in\partial D$ that satisfies \hyperref[eqn:cond1]{\eqref{eqn:cond1}} and \hyperref[eqn:Leb]{\eqref{eqn:Leb}}.
Define $d\tilde{\mu}=|f(\cdot)-f(z)|d\nu+d\mu_{s}$. Then given $\eps>0$ we have
\begin{eqnarray*}
&&\left|\frac{u(x)}{v(x)}-f(z)\right|\\
&=&\left|\frac{1}{v(x)}\left(\int_{\partial D}M_{D}(x,w)(f(w)-f(z))\nu(dw)+\int_{\partial D}M_{D}(x,w)\mu_{s}(dw)\right)\right|\\
&\leq& \frac{\int_{\partial D}M_{D}(x,w)\, \tilde{\mu}(dw)}{v(x)} \\
&=&\frac{\int_{\partial D\cap \{|w-z|\geq \eta \}}M_{D}(x,w) \,\tilde{\mu}(dw)}{v(x)}+\frac{\int_{\partial D}M_{D}(x,w)\,\tilde{\mu}|_{B(z,\eta)}(dw)}{v(x)},
\end{eqnarray*}
where $\tilde{\mu}|_{B(z,\eta)}$ is the truncation of $\tilde{\mu}$ to $B(z,\eta)\cap \partial D$ and $\eta>0$ is a constant which will be determined later. 
Applying Lemma \hyperref[lemma:maximal]{\ref{lemma:maximal}} to the measures $\tilde{\mu}|_{B(z,\eta)}$ and $\nu$, we get
\begin{eqnarray}\label{eqn:RFT}
&&\frac{\int_{\partial D}M_{D}(x,w) \, \tilde{\mu}|_{B(z,\eta)}(dw)}{v(x)}\nonumber\\
&\leq&c_{1}\sup_{r>0}\frac{\tilde{\mu}|_{B(z,\eta)}(B(z,r)\cap \partial D)}{\nu(B(z,r)\cap \partial D)}\nonumber\\
&=&c_{1}\sup_{r\leq \eta}\frac{\int_{\partial D\cap B(z,r)}\left(|f(w)-f(z)|\nu(dw)+\mu_{s}(dw)\right)}{\nu(B(z,r)\cap \partial D)}.
\end{eqnarray}
Using \hyperref[eqn:Leb]{\eqref{eqn:Leb}}, choose $\eta$ so that $\eqref{eqn:RFT}\leq \eps/2$.
Since $|f(\cdot)-f(z)|\in L^{1}(d\nu)$, for this $\eta$ it follows from Lemma \hyperref[lemma:con1]{\ref{lemma:con1}} we can take $\delta$ such that
$$
\left|\frac{\int_{\partial D\cap \{|w-z|\geq \eta \}}M_{D}(x,w) \,\tilde{\mu}(dw)}{v(x)}\right|<\eps/2,
$$
for all $x\in A_{z}^{\beta}$ with $|x-z|<\delta$.
\qed

\section{Harmonic measure and Fatou theorem}\label{Harmonic measure}
In this section we study the harmonic measure that is supported on $\partial D$.
The main result is to show that the harmonic measure supported on $\partial D$ is absolutely continuous with respect to the surface measure of $C^{1,1}$ open sets $D$ and to find the Radon-Nikodym derivative. 

For any Borel subset $A$ of $\R^d$, we use $T_A:=\inf\{t>0: X_t\in A\}$ to denote the first hitting time of $A$.
The next proposition is an analogue of \cite[Proposition 3.1]{Kim}, which was stated only for $x_{0}$ but we remove this restriction and prove the result to hold for all $x\in D$.
\begin{prop}\label{prop:PG}
For any $\lam\in (0, 1/2)$, there exists $c=c(D,d,\phi,\lambda)>0$
such that for any $x,y\in D$ satisfying $|y-x|>2\delta_{D}(y)$ we have
$$
\P_{x}\left(T_{B_{y}^{\lambda}}<\tau_{D}\right)\geq cG_{D}(x,y)\delta_{D}(y)^{d-2},
$$
where $B_{y}^{\lambda}:=B(y,\lambda\delta_{D}(y))$.
\end{prop}
\pf
The proof in the case of $d\ge 3$ is almost identical to that of \cite[Propostion 3.1]{Kim}.
We only give the proof in the case $d=2$.
Since $x\notin B(y, 2\delta_{D}(y))$, $G_{D}(x,\cdot)$
is harmonic in $B(y,2\lambda\delta_{D}(y))$. 
Define $G_{D}1_{B_{y}^{\lambda}}(x):=\int_{B_{y}^{\lambda}}G_{D}(x,z)dz=\E_{x}\left[\int_{0}^{\tau_{D}}1_{B_{y}^{\lambda}}(X_{s})ds\right]$.
By Proposition \hyperref[prop:HP]{\ref{prop:HP}}
there exists a constant $c_{1}>0$ such that
\beq\label{eqn:PG1}
G_{D}1_{B_{y}^{\lambda}}(x)\geq c_{1}G_{D}(x,y)\delta_{D}(y)^{2}.
\eeq
It follows from the strong Markov property that
\beq\label{eqn:PG2}
G_{D}1_{B_{y}^{\lambda}}(x)\leq \P_{x}\left(T_{B_{y}^{\lambda}}
<\tau_{D}\right)\sup_{ w\in \overline{B_{y}^{\lambda}} }
\E_{w}\int_{0}^{\tau_{D}}1_{B_{y}^{\lambda}}\left(X_{s}\right)ds.
\eeq
It follows from \cite[Theorem 1.4]{KSV1} that for any $w\in \overline{B_{y}^{\lambda}}$,
\beq\label{eqn:PG3}
\E_{w}\int_{0}^{\tau_{D}}1_{B_{y}^{\lambda}}\left(X_{s}\right)ds=
\int_{B_{y}^{\lambda}}G_D(w,v)dv\leq c_{2}\int_{B_{y}^{\lambda}}
\ln (1+\frac{\delta_D(w)\delta_D(v)}{|w-v|^2})dv.
\eeq
Note that for $w\in \overline{B_{y}^{\lambda}}$, $\delta_{D}(w)\leq|w-y|+\delta_{D}(y)\leq(1+\lam)\delta_{D}(y)$.
Hence using a polar coordinate system centered at $w$ and integration by parts with
$du=rdr$, $v=\ln(1+\frac{(1+\lam)^{2}\delta_{D}(y)^{2}}{r^{2}})$,
we see that \hyperref[eqn:PG3]{\eqref{eqn:PG3}} is bounded above by
\begin{eqnarray*}
&&c_{2}\int_{B_{y}^{\lambda}}\ln (1+\frac{(1+\lam)^{2}\delta_D(y)^{2}}{|w-v|^2})dv\nonumber\\
&\leq&c_{2}\int_{B(w,2\lam\delta_{D}(y))}\ln (1+\frac{(1+\lam)^{2}\delta_D(y)^{2}}{|w-v|^2})dv\nonumber\\
&\leq&c_{2}\int_{0}^{2\pi}\int_{0}^{2\lam\delta_{D}(y)}r\ln(1+\frac{(1+\lam)^{2}\delta_D(y)^{2}}{r^2})drd\theta\nonumber\\
&\leq&c_{2}(2\pi[\frac{r^{2}}{2}\ln(1+\frac{(1+\lam)^{2}\delta_{D}(y)^{2}}{r^{2}})
]_{0}^{2\lam\delta_{D}(y)}\\
&&+2\pi\int_{0}^{2\lam\delta_{D}(y)}\frac{r(1+\lam)^{2}\delta_{D}(y)^{2}}{r^{2}+(1+\lam)^{2}
\delta_{D}(y)^{2}}dr)\nonumber\\
&\leq&c_{2}(2\pi\frac{(2\lam\delta_{D}(y))^{2}}{2}\ln(1+\frac{(1+\lam)^{2}}{4\lam^{2}})\\
&&+2\pi\left[\frac{1}{2}(1+\lam)^{2}\delta_{D}(y)^{2}\ln(r^{2}+(1+\lam)^{2}\delta_{D}(y)^{2})
\right]_{0}^{2\lam\delta_{D}(y)})\nonumber\\
&\leq&c_{2}\delta_{D}(y)^{2}\left(4\pi\lam^{2}\ln(1+\frac{(1+\lam)^{2}}{4\lam^{2}})+\pi(1+\lam)^{2}
\ln(\frac{4\lam^{2}+(1+\lam)^{2}}{(1+\lam)^{2}})\right)\\
&\leq&c_{3}\delta_{D}(y)^{2}.
\end{eqnarray*}
Combining \hyperref[eqn:PG1]{\eqref{eqn:PG1}}--\hyperref[eqn:PG2]{\eqref{eqn:PG2}} with the display above,
we immediately get the assertion of the proposition.
\qed

For any positive harmonic function $h$ with respect to $X^D$, we use
$\left(\P_{x}^{h},X_{t}^{h}\right)$ to denote the $h$-transform of
$\left(\P_{x},X_{t}^{D}\right)$. That is,
$$
\P_{x}^{h}(A):=\E_{x}\left[\frac{h(X_{t}^{D})}{h(x)};A\right], \quad A\in \mathcal{F}_{t}.
$$
In case $h(\cdot)=M_D(\cdot, z)$ for some $z\in \partial D$, $\left(\P_{x}^{h},X_{t}^{h}\right)$
will be denoted by $\left(\P_{x}^{z},X_{t}^{z}\right)$.
Now we prove a proposition that is an analogue of \cite[Proposition 3.10]{Kim}, which was stated only for $x_{0}$ but we remove this restriction.
\begin{prop}\label{prop:>c}
Suppose that $\lam\in(0,1/2)$.
For any $z\in\partial D$ and $\beta>1$, there exists $c=c(D,d,\phi,\lambda,x,\beta)>0$
such that if $y\in A_{z}^{\beta}$ satisfies $2\delta_D(y)<|x-y|$, then
$$
\P_{x}^{z}\left(T_{B_{y}^{\lambda}}^{z}<\tau_{D}^{z}\right)>c,
$$
where $B_{y}^{\lambda}=B(y,\lambda\delta_{D}(y))$ and $T_{B_{y}^{\lambda}}^{z}:=\inf\{t>0:X_{t}^{z}\in B_{y}^{\lambda}\}$.
\end{prop}
\pf
We only give the proof in the case of $d=2$, the proof in the case $d\ge 3$ is similar.
Fix $z\in\partial D$ and $\beta>1$. Since
$B(y,2\lambda\delta_{D}(y))\subset D$, $M_{D}(\cdot,z)$ is harmonic in $B(y,2\lambda\delta_{D}(y))$.
By the Harnack principle (Proposition \hyperref[prop:HP]{\ref{prop:HP}}), we have
$M_{D}(X_{T_{B_{y}^{\lambda}}},z)\geq c_{1} M_{D}(y,z)$ for some constant $c_{1}>0$.
Now it follows from Proposition \hyperref[prop:PG]{\ref{prop:PG}} that
\begin{eqnarray*}
&&\P_{x}^{z}\left(T^{z}_{B_{y}^{\lam}}<\tau_{D}^{z}\right)\\
&=&\frac{1}{M_{D}(x,z)}\E_{x}\left[M_{D}\left(X_{T_{B^{\lam}_{y}}},z\right),T_{B_{y}^{\lam}}<\tau_{D}\right]\\
&\geq&c_{1}\frac{M_{D}(y,z)}{M_{D}(x,z)}\P_{x}\left(T_{B^{\lam}_{y}}<\tau_{D}\right)\\
&\geq&c_{2}\frac{G_{D}(x,y)M_{D}(y,z)}{M_{D}(x,z)}
\end{eqnarray*}
It follows from \cite[Theorem 1.4]{KSV1} that
$$
G_D(x, y)\ge c_3\ln\left(1+\frac{\delta_D(x)\delta_D(y)}{|x-y|^2}\right).
\quad
$$
Let $\text{diam}D:=\sup\{|x-y| : x,y\in D\}$ be the diameter of a set $D$.
Since $y\in A_{z}^{\beta}$, $|y-z|<\beta \delta_{D}(y)$, $|x-y|\leq \text{diam}D$, and $|x-z|\geq \delta_{D}(x)$. Hence from \hyperref[prop:MK]{\ref{prop:MK}} we have
$$
\P_{x}^{z}\left(T_{B_{y}^{\lambda}}^{z}<\tau_{D}^{z}\right)
\geq c_{4}\frac{\delta_{D}(y)^{2}}{|y-z|^{2}}\frac{|x-z|^2}{|x-y|^2}\geq c_{5}\frac{\delta_{D}(x)^{2}}{\beta^{2}(\text{diam}D)^{2}}.
$$
\qed

Recall that $A\in \mathcal{F}_{\tau_{D}}$ is said to be shift-invariant
if whenever $T<\tau_{D}$ is a stopping time, $1_{A}\circ \theta_{T}=1_{A}$ $\P_{x}$-a.s. for every $x\in D$.
The next proposition is an analogue of \cite[Proposition 3.7]{Kim}.
The proof is identical to that of \cite[Proposition 3.7]{Kim} (see also \cite[p. 196]{Bass}) 
so we omit the proof.

\begin{prop}\label{prop:sh-i}
If $A$ is shift-invariant, then $x\rightarrow \P_{x}^{z}(A)$
is a constant function which is either 0 or 1.
\end{prop}

\begin{prop}\label{prop:exit z}
For any  $z\in \partial D$, we have
$$
\P_{x}^{z}\left(\tau_{D}^{z}<\infty\right)=1, \quad x\in D
$$
and
$$
\P_{x}^{z}\left(\lim_{t\uparrow \tau_{D}^{z}}X_{t}^{z}=z, \tau_{D}^{z}<\infty \right)=1,
\quad x\in D.
$$
\end{prop}
\pf
The proof in the case of $d\ge 3$ is similar to that of \cite[Theorem 3.3]{Kim}.
We only give the proof in the case of $d=2$.
First note that by \cite[Theorem 1.4]{KSV1} and Theorem \hyperref[thm:Martin]{\ref{thm:Martin}} and a similar argument as in \cite[Corollary 6.25]{CZ} we have
$$
\frac{G_{D}(x,y)M_{D}(y,z)}{M_{D}(x,z)}\leq c_1 \left((1\vee \ln(|x-y|^{-1}))+ (1\vee\ln(|y-z|^{-1})\right).
$$
Hence we have
\begin{eqnarray*}
\E_{x}^{z}[\tau_{D}^{z}]&=&\E_{x}^{z}\int_{0}^{\infty}1_{\{t<\tau_{D}^{z}\}}dt\\
&=&\frac{1}{M_{D}(x,z)}\int_{0}^{\infty}\E_{x}\left[M_{D}(X_{t}^{D},z);t<\tau_{D}\right]dt\\
&=&\int_{D}\frac{G_{D}(x,y)M_{D}(y,z)}{M_{D}(x,z)}dy \\
&\leq& c_1\int_{D}\left((1\vee \ln(|x-y|^{-1}))+ (1\vee\ln(|y-z|^{-1})\right)dy<\infty,
\end{eqnarray*}
which implies that $\P_{x}^{z}\left(\tau_{D}^{z}<\infty\right)=1$.

Now we claim that $\P_{x}^{z}\left(\displaystyle\lim_{t\uparrow \tau_{D}^{z}}X_{t}^{z}=z\right)=1$.
Note that the L\'evy process $X$ satisfies the (ACP) condition in \cite[Definition 41.11]{Sa}.
It follows from \cite[Theorem 43.9]{Sa}
that any single point is polar, hence
$\P_y\left(T_{\{x\}}<\infty\right)=0$ for every $x, y\in \R^{d}$.
Now the rest of the proof is the same as that of \cite[Theorem 3.3]{Kim}, \cite[Theorem 5.9]{CZ}, or
\cite[Theorem 3.17]{CS2}.
\qed

The theorem above implies that $\P_{x}^{\cdot}\left(\lim_{t\uparrow\tau_{D}}X_{t}\in K\right)
=1_{K}(\cdot)$ for every $x\in D$ and
Borel subset $K\subset\partial D$.
Hence the next theorem, which is an analogue of \cite[Theorem 3.4]{Kim}, follows easily.
\begin{prop}\label{prop:rep}
Let $\nu$ be a finite measure on $\partial D$. Define
$$
h(x):=\int_{\partial D} M_{D}(x,w)\nu(dw), \quad x\in D.
$$
Then for any $x\in D$ and Borel subset $K$ of $\partial D$,
$$
\P_{x}^{h}\left(\lim_{t\uparrow \tau_{D}^{h}}X_{t}^{h}\in K\right)=\frac{1}{h(x)}\int_{K}M_{D}(x,w)\nu(dw).
$$
\end{prop}

Now the next proposition, which is an analogue of \cite[Proposition 3.5]{Kim},
follows easily from Proposition \hyperref[prop:rep]{\ref{prop:rep}}.
The proof is almost identical to that of  \cite[Proposition 3.5]{Kim} so we omit the proof.

\begin{prop}\label{prop:rep2}
Let $\nu$ be a finite measure on $\partial D$ and $h(x)=\int_{\partial D}M_{D}(x,z)\nu(dz)$.
If $A\in \mathcal{F}_{\tau_{D}}$, then for any Borel subset $K$ of $\partial D$,
$$
\P_{x}^{h}\left(A\cap\{\lim_{t\uparrow\tau_{D}^{h}}X_{t}^{h}\in K\}\right)=\frac{1}{h(x)}\int_{K}\P_{x}^{z}(A)M_{D}(x,z)\nu(dz).
$$
\end{prop}

Now we state a proposition which will play an important role later.
\begin{prop}\label{p:2.12}
Let $u,h$ be nonnegative harmonic functions with respect to $X^{D}$ and $\mu$ and $\nu$ be their Martin measures, respectively.
Let $d\mu=fd\nu +d\mu_{s}$ be Lebesgue decomposition of $\mu$ with respect to
$\nu$.
Then for every $\beta>1$, $x\in D$, and $\nu$-almost every $z\in \partial D$ we have
\beq\label{eqn:2.12}
\P_{x}^{z}\left(\lim_{A^\beta_z\ni x\rightarrow z}\frac{u(x)}{h(x)}=\lim_{t\uparrow \tau_{D}^{z}}\frac{u(X_{t}^{z})}{h(X_{t}^{z})}\right)=1.
\eeq
\end{prop}
\pf
Since $u$ is a nonnegative harmonic function with respect to $X^{D}$,
$u$ is excessive with respect to $X^{D}$. Hence we have
$\E_{x}[u(X^{D}_{t})]\leq u(x)$ for every $x\in D$. So by the Markov property for
the conditioned process, we have for every $t,s>0$
$$
\E_{x}^{h}\left[\frac{u(X_{t+s}^{h})}{h(X_{t+s}^{h})}|\mathcal{F}_{s}\right]
=\E_{X_{s}^{h}}^{h}\left[\frac{u(X_{t}^{h})}{h(X_{t}^{h})}\right]=\frac{1}{h(X_{s}^{h})}\E_{X_{s}^{h}}\left[u(X_{t}^{D})\right]
\leq \frac{u(X_{s}^{h})}{h(X_{s}^{h})}.
$$
Therefore, $u(X_{t}^{h})/h(X_{t}^{h})$ is a nonnegative supermartingale
with respect to $\P_{x}^{h}$ and so by the martingale convergence theorem we have that
$$
\lim_{t\uparrow \tau_{D}^{h}}\frac{u(X_{t}^{h})}{h(X_{t}^{h})}
\text{ exists and is finite } \P_{x}^{h} \text{-a.s.}.
$$
By Proposition \hyperref[prop:rep2]{\ref{prop:rep2}}, we have that
\begin{eqnarray*}
1&=&\P_{x}^{h}\left(\lim_{t\uparrow \tau_{D}^{h}}\frac{u(X_{t}^{h})}{h(X_{t}^{h})}  \text{ exists and is finite }\right)\\
&=&\frac{1}{h(x)}\int_{\partial D}\P_{x}^{z}\left(\lim_{t\uparrow \tau_{D}^{z}}\frac{u(X_{t}^{z})}{h(X_{t}^{z})}  \text{ exists and is finite }\right)M_{D}(x,z)\nu(dz).
\end{eqnarray*}
Since $\P_{x}^{z}\left(\lim_{t\uparrow \tau_{D}^{z}}\frac{u(X_{t}^{z})}{h(X_{t}^{z})}  \text{ exists and is finite }\right)\leq 1$ and $h(x)=\int_{\partial D}M_{D}(x,z)\nu(dz)$, 
we must have
\beq\label{eqn:limit=1}
\P_{x}^{z}\left(\lim_{t\uparrow \tau_{D}^{z}}\frac{u(X_{t}^{z})}{h(X_{t}^{z})}  \text{ exists and is finite }\right)=1,
\eeq
for $\nu-$a.e. $z\in \partial D$.

We will show that \hyperref[eqn:2.12]{\eqref{eqn:2.12}} holds for $z\in\partial D$ satisfying
\hyperref[eqn:cond1]{\eqref{eqn:cond1}}, \hyperref[eqn:Leb]{\eqref{eqn:Leb}},  and \hyperref[eqn:limit=1]{\eqref{eqn:limit=1}}.
For any $\beta>1$,
choose a sequence $y_k\in A^\beta_z$ such that $y_k\to z$.
It follows from Proposition \hyperref[prop:>c]{\ref{prop:>c}} that for any $\lambda\in (0, 1/2)$,
$$
\P_{x}^{z}\left(T_{B_{y_{k}}^{\lambda}}^{z}<\tau_{D}^{z} \text{ i.o.}\right)\geq \liminf_{k\rightarrow \infty}\P_{x}^{z}\left(T_{B_{y_{k}}^{\lambda}}^{z}<\tau_{D}^{z}\right)\geq c>0.
$$
Since $\{T_{B_{y_{k}}^{\lambda}}^{z}<\tau_{D}^{z} \text{ i.o.}\}$ is shift-invariant, by Proposition \hyperref[prop:sh-i]{\ref{prop:sh-i}} we have
$$
\P_{x}^{z}\left(X_{t}^{z} \text{ hits infinitely many } B_{y_{k}}^{\lambda}\right)=\P_{x}^{z}\left(T_{B_{y_{k}}^{\lambda}}^{z}<\tau_{D}^{z} \text{ i.o.}\right)=1.
$$
Suppose that $\{t_{k}, k\in \mathbb{N}\}$ is an increasing sequence
of nonnegative numbers such that $X_{t_{k}}\in B_{y_{k}}^{\lambda}$ under $\P_{x}^{z}$.
By Proposition \hyperref[prop:exit z]{\ref{prop:exit z}} we have $\displaystyle\lim_{t\uparrow \tau_{D}^{z}}X_{t}^{z}=z$ under $\P_{x}^{z}$.
Let $\beta'=(\lambda+\beta)/(1-\lambda)$. Then it is easy to check that
$X_{t_{k}}^{z}\in A^{\beta'}_z$.
Since $\P_{x}^{z}\left(\displaystyle\lim_{k\rightarrow\infty}X_{t_{k}}^{z}=z\right)=1$ it follows from Theorem \hyperref[thm:RFT]{\ref{thm:RFT}}
$$
\P_{x}^{z}\left(\lim_{k\rightarrow \infty}\frac{u(X_{t_{k}}^{z})}{h(X_{t_{k}}^{z})}
=\lim_{A^\beta_z\ni x\rightarrow z}\frac{u(x)}{h(x)}\right)=1.
$$
Since the limit $\displaystyle\lim_{t\uparrow\tau_{D}^{z}}\frac{u(X_{t}^{z})}{h(X_{t}^{z})}$ exists
under $\P_{x}^{z}$, it must be the same as the limit via $t_{k}$. Thus, for any $\beta>1$,
$$
\P_{x}^{z}\left(
\lim_{t\uparrow \tau_{D}^{z}}\frac{u(X_{t}^{z})}{h(X_{t}^{z})}=\lim_{k\rightarrow \infty}\frac{u(X_{t_{k}}^{z})}{h(X_{t_{k}}^{z})}=\lim_{A^\beta_z\ni x\rightarrow z}\frac{u(x)}{h(x)}
\right)=1.
$$
\qed

Now we state the main theorem of this section, which is an analogue of \cite[Thoerem 3.18]{Kim}.
The proof is almost the same with \cite[Thoerem 3.18]{Kim}.
Let $u,h$ be positive harmonic functions with respect to $X^{D}$ and $\mu$ and $\nu$ be their Martin measures, respectively.
Let $d\mu=fd\nu +d\mu_{s}$ be Lebesgue decomposition of $\mu$ with respect to
$\nu$.
Note that it follows from Theorem \hyperref[thm:RFT]{\ref{thm:RFT}} that for any $\beta>1$,
$$
s_{u, h}(z):=\lim_{A_{z}^{\beta}\ni x\rightarrow z}\frac{u(x)}{h(x)}
$$
is well defined for $\nu$-a.e. $z\in \partial D$.

\begin{prop}\label{prop:hrep}
Suppose that $u,h$ are positive harmonic functions with respect to $X^{D}$ and that $u/h$ is bounded.
Let $\nu$ be the Martin measure of $h$.
For every $x\in D$ we have
$$
u(x)=\int_{\partial D}M_{D}(x,z)s_{u, h}(z)\nu(dz).
$$
Equivalently, $s_{u, h}(z)$ is the Radon-Nikodym derivative of
the Martin measure of $u$ with respect to $\nu$.
\end{prop}
\pf
It follows from Proposition \hyperref[p:2.12]{\ref{p:2.12}} that for every $x\in D$ and $\nu$-a.e. $z\in\partial D$ and $\beta>1$,
$$
\P_{x}^{z}\left(\lim_{A_{z}^{\beta}\ni x\rightarrow z}\frac{u(x)}{h(x)}=\lim_{t\uparrow \tau_{D}^{z}}\frac{u(X_{t}^{z})}{h(X_{t}^{z})}\right)=1.
$$
Now take an increasing sequence of smooth open sets $\{D_{n}\}_{n\geq 1}$ such that
$\overline{D_{n}}\subset D_{n+1}$ and $\cup_{n=1}^{\infty}D_{n}=D$. Then we have
\begin{eqnarray*}
1&=&\P_{x}^{z}\left(\lim_{n\rightarrow \infty}\left(\frac{u}{h}\right)\left(X_{\tau_{n}^{z}}^{z}\right)
=\lim_{t\uparrow\tau_{D}^{z}}\frac{u(X_{t}^{z})}{h(X_{t}^{z})}
=\lim_{A_{z}^{\beta}\ni x\rightarrow z}
\frac{u(x)}{h(x)}\right)\\
&=&\P_{x}^{z}\left(\lim_{n\rightarrow \infty}\left(\frac{u}{h}\right)\left(X_{\tau_{n}^{z}}^{z}\right)=s_{u, h}(z),
\lim_{t\uparrow \tau_{D}^{z}}X_{t}^{z} =z\right)\\
&=&\P_{x}^{z}\left(\lim_{n\rightarrow \infty}\left(\frac{u}{h}\right)\left(X_{\tau_{n}^{z}}^{z}\right)
=s_{u, h}\left(\lim_{t\uparrow \tau_{D}^{z}}X_{t}^{z}\right)\right)
\end{eqnarray*}
for $\nu$-a.e. $z\in\partial D$.
By Propositions  \hyperref[prop:sh-i]{\ref{prop:sh-i}} and \hyperref[prop:rep2]{\ref{prop:rep2}} we have
\begin{eqnarray*}
1&=&\frac{1}{h(x)}\int_{\partial D}\P_{x}^{z}\left(\lim_{n\rightarrow \infty}\left(\frac{u}{h}\right)\left(X_{\tau_{n}^{z}}^{z}\right)
=s_{u, h}\left(\lim_{t\uparrow \tau_{D}^{z}}X_{t}^{z}\right)\right)M_{D}(x,z)\nu(dz)\\
&=&\P_{x}^{h}\left(\lim_{n\rightarrow\infty}\left(\frac{u}{h}\right)\left(X_{\tau_{D_{n}}^{h}}^{h}\right)=
s_{u, h}\left(\lim_{t\uparrow\tau_{D}^{h}}X_{t}^{h}\right)\right).
\end{eqnarray*}
Therefore, by the bounded convergence theorem and the harmonicity of $u/h$ with respect to $\P_{x}^{h}$, we have
$$
\frac{u(x)}{h(x)}=\lim_{n\rightarrow \infty}\E_{x}^{h}\left[\left(\frac{u}{h}\right)\left(X_{\tau_{n}^{h}}^{h}\right)\right]
=\E_{x}^{h}\left[\lim_{n\rightarrow \infty}\left(\frac{u}{h}\right)\left(X_{\tau_{n}^{h}}^{h}\right)\right]
=\E_{x}^{h}\left[s_{u, h}\left(\lim_{t\uparrow\tau_{D}^{h}}X_{t}^{h}\right)\right]
$$
for every $x\in D$.
By Proposition \hyperref[prop:rep]{\ref{prop:rep}} we have
\beq\label{eqn:hrep1}
\E_{x}^{h}\left[1_{K}\left(\lim_{t\uparrow\tau_{D}^{h}}X_{t}^{h}\right)\right]=\frac{1}{h(x)}\int_{\partial D}M_{D}(x,w)1_{K}(w)\nu(dw).
\eeq
Clearly \hyperref[eqn:hrep1]{\eqref{eqn:hrep1}} remains true if $1_{K}(w)$ is replaced by simple functions of the form $\sum_{i=1}^{n}a_{i}1_{A_{i}}(w)$ where $a_{i}\geq 0$
and $A_{i}\subset\partial D$ are disjoint Borel subsets of $\partial D$.
Since $s_{u, h}$ is bounded, there exists a sequence of bounded simple functions
$f_{n}(w)\leq s_{u, h}(w)$ converging to $s_{u, h}(w)$. Then it follows from bounded convergence theorem that
\begin{eqnarray*}
&&\frac{u(x)}{h(x)}=\E_{x}^{h}\left[s_{u,h}\left(\lim_{t\uparrow\tau_{D}^{h}}X_{t}^{h}\right)\right]=\lim_{n\rightarrow \infty}\E_{x}^{h}\left[f_{n}\left(\lim_{t\uparrow\tau_{D}^{h}}X_{t}^{h}\right)\right]\\
&=&\lim_{n\rightarrow \infty}\frac{1}{h(x)}\int_{\partial D}M_{D}(x,w)f_{n}(w)\nu(dw)\\
&=&\frac{1}{h(x)}\int_{\partial D}M_{D}(x,w)s_{u, h}(w)\nu(dw).
\end{eqnarray*}
Now the proof is complete.
\qed

In order to study the harmonic measure supported on $\partial D$, we need auxiliary functions. 
Let
$$
F(x):=\P_{x}\left(X_{\tau_{D}}\in\partial D\right), \quad x\in D
$$
and
$$
G(x)=\int_{\partial D}M_{D}(x,z)\sigma(dz), \quad x\in D
$$
where $\sigma$ is the surface measure of $\partial D$.
It is easy to see that $F(x)$ and $G(x)$ are harmonic with respect to $X^{D}$.
Now we prove that $G(x)$ is bounded on $D$.

\begin{lemma}\label{lemma:bdd}
There exist constants $C_{3}, C_{4}$ depending only on $D, d, \phi , x_{0}$ such that
$$
0<C_{3}\leq G(x) \leq C_{4} <\infty.
$$
\end{lemma}
\pf
Recall that $D$ satisfies the Ahlfors regular condition \hyperref[eqn:Ahlfors]{\eqref{eqn:Ahlfors}}.
First suppose that $\delta_{D}(x)\geq R_{1}$. Then we have $\text{diam} D\geq |x-z|\geq \delta_{D}(x)\geq R_{1}$ for any $z\in\partial D$. 
Hence it follows from Proposition \hyperref[prop:MK]{\ref{prop:MK}} we have
$M_{D}(x,z)\geq c_{1}\frac{\delta_{D}(x)}{|x-z|^{d}}\geq c_{1}\frac{R_{1}}{(\text{diam} D)^{d}}$ and 
$M_{D}(x,z)\leq c_{2}\frac{\delta_{D}(x)}{|x-z|^{d}}\leq c_{2}\delta_{D}(x)^{1-d}\leq c_{2}R_{1}^{1-d}$.
Hence we have
$$
G(x)=\int_{\partial D}M_{D}(x,z)\sigma(dz)\geq c_{1}\frac{R_{1}}{(\text{diam} D)^{d}}\sigma(\partial D), 
$$
and
$$
G(x)=\int_{\partial D}M_{D}(x,z)\sigma(dz)\leq c_{2}\sigma(\partial D)R_{1}^{1-d}.
$$

Now suppose that $\delta_{D}(x)<R_{1}$.
For each $x\in D$ let $P=P(x)\in \partial D$ be a point such that $|x-P|=\delta_{D}(x)$.
Let $A_{n}=A_{n}(x)=\{z\in\partial D : 2^{n-1}|x-P|\leq |x-z| < 2^{n} |x-P|\}$, $n\in \mathbb{N}$.
Since $D$ is bounded, there exists $N=N(x)$
such that $\partial D\subset \bigcup_{n=1}^{N}A_{n}$. 
Note that $\{z\in\partial D: |z-P|<|x-P|\}\subset A_{1}$ since if $|z-P|<|x-P|$ then $|x-z|\leq |x-P|+|P-z|< 2|x-P|$ and $|x-z|\geq |x-P|$ for any $z\in\partial D$.
Since $\delta_{D}(x)=|x-P|<R_{1}$ it follows from \hyperref[eqn:Ahlfors]{\eqref{eqn:Ahlfors}}
\begin{eqnarray*}
G(x)&=&\int_{\partial D}M_{D}(x,z)\sigma(dz)\geq\int_{A_1} M_{D}(x,z)\sigma(dz)\\
&\geq&c_{3}\int_{A_1}\frac{\delta_{D}(x)}{|x-z|^{d}}\sigma(dz)\,\geq\, c_{4}\int_{A_1}\frac{\delta_{D}(x)}{2^{d}|x-P|^{d}}\sigma(dz)\\
&\geq&c_{4}\int_{\{z\in\partial D: |z-P|<|x-P|\}}\frac{\delta_{D}(x)}{2^{d}|x-P|^{d}}\sigma(dz)\\
&\geq&c_{4}\frac{\delta_{D}(x)}{2^{d}|x-P|^{d}}\sigma(\{z\in\partial D : |z-P|<|x-P|\})\\
&\geq&c_{5}\frac{\delta_{D}(x)}{2^{d}|x-P|^{d}}|x-P|^{d-1}\\
&\geq&\frac{c_{5}}{2^d}.
\end{eqnarray*}

Now we prove an upper bound. Notice that for any $0<r<\text{diam} D$ there exists a constant $c_{6}$ such that $\sigma(\partial D\cap B(z,r))\leq c_{6}r^{d-1}$.
If $r<R_{1}$ this is just \hyperref[eqn:Ahlfors]{\eqref{eqn:Ahlfors}}. If $r\geq R_{1}$ then $\sigma(\partial D\cap B(z,r))\leq \sigma(\partial D)= c_{6}R_{1}^{d-1}\leq c_{6}r^{d-1}$, where $c_{6}:=\frac{\sigma(\partial D)}{R_{1}^{d-1}}$. Since $A_{n}\subset \{z\in\partial D : |z-P|\leq (2^{n}+1)|x-P|\}\subset \{z\in\partial D : |z-P|\leq 2^{n+1}|x-P|\}$. 
\begin{eqnarray*}
G(x)&=&\int_{\partial D}M_{D}(x,z)\sigma(dz)\,\leq\,\sum_{n=1}^{N}\int_{A_{n}}M_{D}(x,z)\sigma(dz)\\
&\leq&c_{7}\sum_{n=1}^{N}\int_{A_{n}}\frac{\delta_{D}(x)}{|x-z|^{d}}\sigma(dz)\,\leq\, c_{7}\sum_{n=1}^{N}\int_{A_{n}}\delta_{D}(x)\left(2^{n-1}|x-P|\right)^{-d}\sigma(dz)\\
&\leq&c_{7}\,\delta_{D}(x)^{1-d}\sum_{n=1}^{N}2^{-d(n-1)}\sigma( A_n )\\
&\leq&c_{7}\delta_{D}(x)^{1-d}\sum_{n=1}^{N}2^{-d(n-1)}\sigma\left(\{z\in\partial D : |z-P|\leq 2^{n+1}|x-P|\}\right)\\
&\leq&c_{8}\delta_{D}(x)^{1-d}\sum_{n=1}^{N}2^{-d(n-1)}(2^{n+1}|x-P|)^{d-1}\\
&\leq&c_{8}2^{2d-1}\sum_{n=1}^{N}2^{-n}\leq c_{8}2^{2d-1}\sum_{n=1}^{\infty}2^{-n}=c_{8}2^{2d-1}.
\end{eqnarray*}
Now set $C_{3}:=c_{1}\frac{R_{1}}{(\text{diam}D)^{d}}\sigma(\partial D)\wedge \frac{c_{5}}{2^{d}}$ and $C_{4}:=c_{2}R_{1}^{1-d}\vee c_{8}2^{2d-1}$.
\qed

It follows from Lemma \hyperref[lemma:bdd]{\ref{lemma:bdd}} that $\frac{F(x)}{G(x)}$ is bounded in $D$.
Thus it follows from Proposition \hyperref[prop:hrep]{\ref{prop:hrep}} that, for any $\beta>1$, the limit
\begin{equation}\label{eqn:defns}
s_{F, G}(z)=\lim_{A^\beta_z\ni x\rightarrow z}\frac{F(x)}{G(x)}
\end{equation}
exists $\sigma -$a.e. $z\in\partial D$ and $F$ can be written as
$$
F(x)=\int_{\partial D}M_{D}(x,w)s_{F,G}(w)\sigma(dw), \quad x\in D.
$$

Next proposition says as the starting point $x$ approaches $\partial D$, the probability that subordinate Brownian motions with Gaussian components exit the domain through the boundary of the domain $\partial D$ converges to 1. It was proved in a more general setting in \cite[Theorem 3.2]{Mi} and we record the fact here for the reader's convenience. 

\begin{prop}\label{prop:limit=1}
Let $D$ be a $C^{1,1}$ open set in $\R^{d}$, $d\geq 2$. Then for every point $z\in\partial D$
$$
\lim_{D\ni x\rightarrow z\in\partial D}\P_{x}\left(X_{\tau_{D}}\in\partial D\right)=1.
$$
\end{prop}
\pf
Subordinate Brownian motions are isotropic processes hence they satisfy conditions $(H_{1};\R^{d},\alpha)$ and $(H_{2};\R^{d},\alpha)$ in \cite{Mi}
and all points of $D$ are possible (see \cite{Mi} for details). Hence it follows from the remark $(c)$ under \cite[Theorem 3.2]{Mi} $\displaystyle\lim_{x\rightarrow z\in \partial D}\P_{x}\left(X_{\tau_{D}}\in\partial D\right)=1$.  
\qed

Now the next result follows immediately from Proposition \hyperref[prop:hrep]{\ref{prop:hrep}}.
\begin{thm}\label{thm:Rep harmonic boundary}
For any $\beta>1$, the limit
$$
s_{F, G}(z)=\displaystyle\lim_{A^\beta_z\ni x\rightarrow z\in \partial D}\frac{F(x)}{G(x)}
=\lim_{A^\beta_z\ni x\rightarrow z\in \partial D}\frac{1}{G(x)}
$$
exists $\sigma -$a.e. $z\in\partial D$ and
$0<C_{3}\leq s_{F,G}(z)\leq C_{4}<\infty$.
Furthermore, $F(x)$ can be written as
$$
F(x)=\int_{\partial D}M_{D}(x,w)s_{F,G}(w)\sigma(dw), \quad x\in D.
$$
\end{thm}

As a corollary of
Proposition \hyperref[prop:limit=1]{\ref{prop:limit=1}}, Theorem \hyperref[thm:Rep harmonic boundary]{\ref{thm:Rep harmonic boundary}} and \hyperref[thm:RFT]{\ref{thm:RFT}} we can prove Fatou's theorem for nonnegative harmonic functions with respect to $X^{D}$.
\begin{corollary}\label{cor:Fatou killed}
Let $u(x)$ be nonnegative and harmonic with respect to $X^{D}$ on $D$. Then for any $\beta>1$ the nontangential limit
$$
\lim_{A_{z}^{\beta}\ni x\rightarrow z}u(x)
$$
exists for $\sigma-$a.e. $z\in \partial D$.
\end{corollary}
\pf
From Theorem \hyperref[thm:Rep harmonic boundary]{\ref{thm:Rep harmonic boundary}} we have $F(x)=\int_{\partial D}M_{D}(x,w)s_{F,G}(w)\sigma(dw)$, $x\in D$. 
It follows from Theorems \hyperref[thm:RFT]{\ref{thm:RFT}} and \hyperref[thm:Rep harmonic boundary]{\ref{thm:Rep harmonic boundary}} $\displaystyle\lim_{A_{z}^{\beta}\ni x\rightarrow z}\frac{u(x)}{F(x)}$ exists $\sigma$-a.e. $z\in\partial D$.
From Proposition \hyperref[prop:limit=1]{\ref{prop:limit=1}}  we have $\displaystyle\lim_{x\rightarrow z}F(x)=1$. Hence 
$$
\lim_{A_{z}^{\beta}\ni x\rightarrow z}u(x) \text{ exists } \sigma\text{-a.e. } z\in\partial D.
$$
\qed

Now we show that
\begin{equation}\label{eqn:pk}
P_D(x, z):=M_{D}(x,z)s_{F,G}(z), \quad z\in \partial D
\end{equation}
is the Radon-Nikodym derivative of the restriction of harmonic measure $\P_{x}(X_{\tau_{D}}\in \cdot)$ to $\partial D$ with respect to the
surface measure $\sigma$ on $\partial D$.
In order to do this, we need a few lemmas.
For any $z\in\partial D$, we let $\phi_{z}$ be the $C^{1,1}$ function associated with $z$ in the definition
of $C^{1,1}$ open set. For any $x\in \{y=(\tilde{y},y_{d})\in B(z,R_{0}): y_{d}>\phi_{z}(\tilde{y})\}$ we put
$\rho_{z}(x):=x_{d}-\phi_{z}(\tilde{x})$.
For $r_{1},r_{2}>0$, we define
$$
D_{z}(r_{1},r_{2}):=\{y\in D : r_{1}>\rho_{z}(y)>0, |\tilde{y}|<r_{2}\}.
$$
Let $R_{2}:=R_{0}/4(\sqrt{1+(1+\Lambda_{0})^{2}})$.
The following result is \cite[Lemma 4.3]{KSV1}.
\begin{lemma}\label{lemma:exit2}\emph{(\cite[Lemma 4.3]{KSV1})}
There exist constants $\lambda_{0}>2R_{2}^{-1}$, $\kappa_{0}\in(0,1)$
and $c=c(R_{0},\Lambda_{0})$ such that for every
$\lam\geq\lambda_{0}$, $z\in\partial D$ and $x\in D_{z}(2^{-1}(1+\Lambda_{0})^{-1}
\kappa_{0}\lam^{-1},\kappa_{0}\lam^{-1})$
with $\tilde{x}=0$,
$$
\P_{x}\left(X_{\tau_{D_{z}(\kappa_{0}\lambda^{-1},\lambda^{-1})}}\in D\right)\leq c\lambda\delta_{D}(x).
$$
\end{lemma}
\begin{lemma}\label{lemma:exit}
For any $r<R_{0}$ and $z\in\partial D$ we have
$$
\lim_{x\rightarrow z}\P_{x}\left(X_{\tau_{D}}\notin B(z,r), X_{\tau_{D}}\in \partial D\right)=0.
$$
That is, for any $\eps>0$ there exists a constant $\delta=\delta(\eps)>0$
such that for any $x\in D$ with $|x-z|<\delta$
$
\P_{x}\left(X_{\tau_{D}}\notin B(z,r), X_{\tau_{D}}\in \partial D\right)<\eps.$
\end{lemma}
\pf
For any $r<R_{0}$, we take a large enough $\lambda$ so that
$D_{z}(\kappa_{0}\lambda^{-1},\lambda^{-1})\subset B(z,r)$.
Then,
$$
\{X_{\tau_{D}}\notin B(z,r), X_{\tau_{D}}\in \partial D\}\subset \{X_{\tau_{D_{z}(\kappa_{0}\lambda^{-1},\lambda^{-1})}} \in D\}.
$$
It follows from Lemma \hyperref[lemma:exit2]{\ref{lemma:exit2}} that
\begin{eqnarray*}
&&\P_{x}\left(X_{\tau_{D}}\notin B(z,r), X_{\tau_{D}}\in \partial D\right)\leq \P_{x}\left(X_{\tau_{D_{z}(\kappa_{0}\lambda^{-1},\lambda^{-1})}} \in D\right)\\
&\leq& c\lambda \delta_{D}(x)\leq c\lambda |x-z|.
\end{eqnarray*}
By taking $\delta=\eps/c\lambda$, we arrive at the desired assertion.
\qed
\begin{lemma}\label{lemma:conv cts}
For any continuous function $g$ on $\partial D$, define
$$
u_{g}(x):=\E_{x}\left[g(X_{\tau_{D}}), X_{\tau_{D}}\in \partial D\right], \quad x\in D.
$$
Then for any $z\in\partial D$,
$$
\lim_{x\rightarrow z\in \partial D}u_{g}(x)=g(z).
$$
Furthermore $u_{g}(x)$ is given by
$$
u_{g}(x)=\int_{\partial D}M_{D}(x,w)s_{F,G}(w)g(w)\sigma(dw), \quad x\in D.
$$
\end{lemma}
\pf
For any $\eps>0$, let $\delta_{1}=\delta_{1}(\eps)>0$ be such that
\beq\label{eqn:exit2}
|g(y)-g(z)|<\eps \text{ whenever } |y-z|\leq\delta_{1}.
\eeq
Without loss of generality we may assume $\delta_{1}<R_{0}$.
Let $\delta$ be the constant in Lemma \hyperref[lemma:exit]{\ref{lemma:exit}} so that
\beq\label{eqn:exit3}
\P_{x}\left(X_{\tau_{D}}\notin B(z,\delta_{1}), X_{\tau_{D}}\in \partial D\right)<\eps
\eeq
for $|x-z|<\delta$.
It follows from Proposition \hyperref[prop:limit=1]{\ref{prop:limit=1}} that there exists $\delta_{2}>0$ such that
\beq\label{eqn:exit4}
\P_{x}\left(X_{\tau_{D}}\notin\partial D\right)<\eps
\eeq
for $|x-z|<\delta_{2}$.
Combining \hyperref[eqn:exit2]{\eqref{eqn:exit2}}--\hyperref[eqn:exit4]{\eqref{eqn:exit4}} we get that,
for any $x$ satisfying $|x-z|<\delta\wedge \delta_{2}$,
\begin{eqnarray*}
&&\left|u_{g}(x)-g(z)\right|\\
&=&\left|\E_{x}\left[g(X_{\tau_{D}}), X_{\tau_{D}}\in \partial D\right]-g(z)\right|\\
&=&|\E_{x}\left[g(X_{\tau_{D}}), X_{\tau_{D}}\in B(z,\delta_{1}), X_{\tau_{D}}\in \partial D\right]\\
&&+ \E_{x}\left[g(X_{\tau_{D}}), X_{\tau_{D}}\notin B(z,\delta_{1}),X_{\tau_{D}}\in \partial D\right]\\
&&-g(z)\P_{x}\left(X_{\tau_{D}}\in B(z,\delta_{1}), X_{\tau_{D}}\in \partial D\right)\\
&&-g(z)\P_{x}\left(X_{\tau_{D}}\notin B(z,\delta_{1}), X_{\tau_{D}}\in \partial D\right)-g(z)\P_{x}\left(X_{\tau_{D}}\notin \partial D\right)|\\
&\leq& 2\|g\|_{\infty}\P_{x}\left(X_{\tau_{D}}\notin B(z,\delta_{1}), X_{\tau_{D}}\in \partial D\right) +\|g\|_{\infty}\P_{x}\left( X_{\tau_{D}}\notin \partial D\right)\\
&&+\E_{x}\left[\left|g(X_{\tau_{D}})-g(z)\right|,X_{\tau_{D}}\in \partial D, X_{\tau_{D}}\in B(z,\delta_{1})\right]\\
&\leq&4\eps\|g\|_{\infty}.
\end{eqnarray*}

It follows from \hyperref[eqn:defns]{\eqref{eqn:defns}} $\displaystyle\lim_{x\rightarrow z\in\partial D}\frac{u_{g}(x)}{G(x)}=g(z)s_{F,G}(z)$ for $\sigma$-a.e. $z\in\partial D$. 
Since $D$ is bounded $\partial D$ is compact and $u_{g}(x)$ is bounded. 
Hence from Proposition \hyperref[prop:hrep]{\ref{prop:hrep}} we have
$$
u_{g}(x)=\int_{\partial D}M_{D}(x,w)s_{F,G}(w)g(w)\sigma(dw), \quad x\in D.
$$
\qed

\begin{thm}\label{thm:har rep}
For any (Lebesgue) measurable set $A\in \partial D$
$$
\P_{x}\left(X_{\tau_{D}}\in A\right)=\int_{A}M_{D}(x,w)s_{F,G}(w)\sigma(dw), \quad x\in D.
$$
\end{thm}
\pf
Let $A$ be a (Lebesgue) measurable set in $\partial D$.
Choose bounded and continuous functions $f_{n}(x)$ converging to $1_{A}(x)$.
Then it follows from Lemma \hyperref[lemma:conv cts]{\ref{lemma:conv cts}} and the dominated convergence theorem
\begin{eqnarray*}
&&\P_{x}(X_{\tau_{D}}\in A)=\E_{x}(1_{A}(X_{\tau_{D}}))=\E_{x}(\lim_{n\rightarrow\infty}f_{n}(X_{\tau_{D}}))=\lim_{n\rightarrow\infty}\E_{x}(f_{n}(X_{\tau_{D}}))\\
&=&\lim_{n\rightarrow\infty}\int_{\partial D}M_{D}(x,w)s_{F,G}(w)f_{n}(w)\sigma(dw)\\
&=&\int_{\partial D}M_{D}(x,w)s_{F,G}(w)1_{A}(w)\sigma(dw).
\end{eqnarray*}
\qed

Combining this with the two-sided estimates on $M_D(x,z)$, we have the following theorem.

\begin{thm}\label{thm:sharp Poisson}
The function $P_D$ defined in \hyperref[eqn:pk]{\eqref{eqn:pk}} is
the Radon-Nikodym derivative of the restriction of the harmonic measure $\P_{x}(X_{\tau_{D}}\in \cdot )$ to $\partial D$ with respect to the
surface measure $\sigma$ on $\partial D$. Furthermore, there exist positive constants $C_{3}(D,d,\phi,x_{0})<C_{4}(D,d,\phi,x_{0})$
such that
$$
C_3\frac{\delta_D(x)}{|x-z|^d}\le P_D(x, z)\le C_4\frac{\delta_D(x)}{|x-z|^d},
\quad (x, z)\in D\times \partial D.
$$
Therefore the harmonic measure restricted to $\partial D$ is mutually absolutely continuous with respect to the surface measure $\sigma$ on $\partial D$.
\end{thm}

\section{Integral representation of harmonic functions with respect to $X$}\label{Main result}
In this section we investigate the integral representation of nonnegative harmonic functions with respect to $X$ and show that tangential convergence of harmonic functions with respect to $X^{D}$ can fail.

Let $D$ be a bounded $C^{1,1}$ open set in $\R^{d}$, $d\geq 2$. Take a sequence of smooth open sets $\{D_{n}\}$, $D_{n}\subset \overline{D_{n}}\subset D_{n+1}$, $\cup_{n=1}^{\infty}D_{n}=D$. Let $\tau_{n}:=\tau_{D_{n}}$ be the first exit time of $D_{n}$.
We will define some auxiliary sets $\mathcal{A},\mathcal{B},\mathcal{C},\mathcal{D}\subset \Omega$. Let
$$
\mathcal{A}=\{w\in \Omega : X_{\tau_{D}^{-}}\neq X_{\tau_{D}}\},  \mathcal{B}=\Omega\setminus \mathcal{A}, 
$$
$$
\mathcal{C}=\{\omega\in \Omega : \tau_{n}=\tau_{D}
\text{ for some }n\in \N\}, \text{and } \mathcal{D}=\Omega\setminus \mathcal{C}.
$$
Since $X_{\tau_{D}}\in D^{c}$ and $X_{\tau_{D}^{-}}\in \overline{D}$ for $x\in D$ we have $\P_{x}$-almost surely 
$$
\mathcal{A}=\{\omega\in\Omega : X_{\tau_{D}} \in \overline{D}^{c}\}=\{w\in \Omega : X_{\tau_{D}}\notin\partial D\}.
$$
Suppose that $\omega \in \mathcal{A}\setminus \mathcal{C}$. 
Then $\tau_{n}(\omega)<\tau_{D}(\omega)$ for all $n\in \N$. By the quasi-left continuity of L\'evy processes we have $\displaystyle\lim_{\tau_{n}\uparrow \tau_{D}}X_{\tau_{n}}(\omega)=X_{\tau_{D}}(\omega)$. But this implies $X_{\tau_{D}}(\omega)\in D^{c}\cap \overline{D}=\partial D$, which is a contradiction. Hence $\mathcal{A}\setminus \mathcal{C}=\emptyset$ or
$\mathcal{A}\subset \mathcal{C}$. By taking complement we also have $\mathcal{D}\subset \mathcal{B}$. 

Finally consider $\mathcal{C}\setminus \mathcal{A}=\{\omega\in\Omega : \tau_{n}=\tau_{D} \text{ for some } n \text{ and } X_{\tau_{D}}\in \partial D\}$.
Clearly $\mathcal{C}\setminus \mathcal{A}\subset \{X_{\tau_{n}}\in \partial D \text{ for some }n\}$.
Note that $\{X_{\tau_{n}}\in \partial D\}=\{X_{\tau_{n}}\in \partial D \text{ and } X_{\tau_{n}^{-}}\neq X_{\tau_{n}}\}$ since $\overline{D}_{n}\subset D$.
Since $|\partial D|=0$ for any $C^{1,1}$ open set $D$, it follows from \hyperref[eqn:Levy system2]{\eqref{eqn:Levy system2}} we have
$$
\P_{x}(X_{\tau_{n}}\in \partial D)=\P_{x}(X_{\tau_{n}}\in \partial D, X_{\tau_{n}^{-}}\neq X_{\tau_{n}})=\int_{\partial D}K_{D_{n}}(x,z)dz=0.
$$
Hence we conclude that under $\P_{x},$ $x\in D$
$$
\P_{x}(\mathcal{C}\setminus \mathcal{A})=0.
$$
Hence from now on we will identify all these sets to be equal under $\P_{x}$. That is we let 
$$
\{w\in \Omega : X_{\tau_{D}^{-}}\neq X_{\tau_{D}}\}=\{\omega\in\Omega : X_{\tau_{D}}\in \overline{D}^{c}\}=\{\omega\in\Omega :  \tau_{n}=\tau_{D} \text{ for some } n\},
$$
$$
\{w\in \Omega : X_{\tau_{D}^{-}}=X_{\tau_{D}}\}=\{\omega\in\Omega : X_{\tau_{D}}\in  \partial D\}=\{\omega\in\Omega :  \tau_{n}\neq \tau_{D} \text{ for all } n\}.
$$

Let $\mathcal{J}$ be
$$
\mathcal{J}=\{\tau_{n}=\tau_{D} \text{ for some } n\}
=\{ X_{\tau_{D}} \in \overline{D}^{c}\}
=\{X_{\tau_{D}^{-}}\neq X_{\tau_{D}}\}.
$$

\begin{lemma}\label{lemma:rep jump2}
Let $u$ be a nonnegative function defined on $D$. For any open set $B\subset \overline{B}\subset D$ we have
$$
\E_{x}[u(X_{\tau_{D}}), \mathcal{J}]=\E_{x}[u(X_{\tau_{B}}), \tau_{B}=\tau_{D}]+\E_{x}[\E_{X_{\tau_{B}}}\left(u(X_{\tau_{D}}), \mathcal{J}\right),\tau_{B}<\tau_{D}].
$$
\end{lemma}
\pf
Take an increasing sequence of smooth opens sets $\{D_{n}\}$ as in the beginning of the chapter.
For any open set $B\subset \overline{B}\subset D$ we can take an open set $D_{k}$ such that $B\subset D_{k}$. Then $\tau_{B}\leq \tau_{D_{k}}$. Hence we have
$\{\tau_{B}=\tau_{D}\}\subset \{\tau_{D_{k}}=\tau_{D}\} \subset \mathcal{J}$. Hence
it suffices to show
$$
\E_{x}[u(X_{\tau_{D}}), \mathcal{J}\setminus \{\tau_{B}=\tau_{D}\}]
=\E_{x}[\E_{X_{\tau_{B}}}\left(u(X_{\tau_{D}}), \mathcal{J}\right),\tau_{B}<\tau_{D}].
$$
From the strong Markov property of $X$ we have
\begin{eqnarray*}
&&\E_{x}[u(X_{\tau_{D}}), \mathcal{J}\setminus \{\tau_{B}=\tau_{D}\}]\\
&=&\E_{x}[u(X_{\tau_{D}}), \mathcal{J}\cap \{\tau_{B}<\tau_{D}\}]\\
&=&\E_{x}[\E[u(X_{\tau_{D}}), \mathcal{J}\cap \{\tau_{B}<\tau_{D}\}| \mathcal{F}_{\tau_{B}}]]\\
&=&\E_{x}[\E[u(X_{\tau_{D}}), \mathcal{J} | \mathcal{F}_{\tau_{B}}], \tau_{B}<\tau_{D}]\\
&=&\E_{x}[\E_{X_{\tau_{B}}}\left(u(X_{\tau_{D}}),\mathcal{J}\right),\tau_{B}<\tau_{D}].
\end{eqnarray*}
\qed

\begin{lemma}\label{lemma:rep jump3}
Let $u$ be a harmonic function on $D$ with respect to $X$. Let $v(x):=u(x)-\E_{x}[u(X_{\tau_{D}}), \mathcal{J}]$.
Then $v(x)$ is nonnegative and harmonic with respect to $X^{D}$ on $D$.
\end{lemma}
\pf
Take $D_{n}\subset \overline{D_{n}}\subset D_{n+1}\uparrow D$. Since $u$ is harmonic with respect to $X$ we have
$$
u(x)=\E_{x}[u(X_{\tau_{n}})]\geq \E_{x}[u(X_{\tau_{D}}), \tau_{n}=\tau_{D}].
$$
As $n\rightarrow \infty$ we have $\{\tau_{n}=\tau_{D}\}\uparrow \mathcal{J}$. By the monotone convergence theorem we have
$$
u(x)\geq \E_{x}[u(X_{\tau_{D}}), \mathcal{J}].
$$

From the harmonicity of $u$ and Lemma \hyperref[lemma:rep jump2]{\ref{lemma:rep jump2}}, for any open set $B\subset D$ whose closure is compact in $D$ we have
\begin{eqnarray*}
&&\E_{x}[v(X_{\tau_{B}}^{D})]\\
&=&\E_{x}[u(X_{\tau_{B}}^{D})]-\E_{x}\left[\E_{X_{\tau_{B}}^{D}}[u(X_{\tau_{D}}), \mathcal{J}]\right]\\
&=&\E_{x}[u(X_{\tau_{B}}), \tau_{B}<\tau_{D}]-\E_{x}\left[\E_{X_{\tau_{B}}}[u(X_{\tau_{D}}), \mathcal{J}], \tau_{B}<\tau_{D}\right]\\
&=&\E_{x}[u(X_{\tau_{B}})]-\E_{x}[u(X_{\tau_{B}}), \tau_{B}=\tau_{D}]-\E_{x}\left[\E_{X_{\tau_{B}}}[u(X_{\tau_{D}}), \mathcal{J}], \tau_{B}<\tau_{D}\right]\\
&=&u(x)-\E_{x}[u(X_{\tau_{D}}), \mathcal{J}]+\E_{x}[u(X_{\tau_{D}}), \mathcal{J}]\\
&&-\E_{x}[u(X_{\tau_{B}}), \tau_{B}=\tau_{D}]-\E_{x}\left[\E_{X_{\tau_{B}}}[u(X_{\tau_{D}}), \mathcal{J}], \tau_{B}<\tau_{D}\right]\\
&=&v(x)+\E_{x}[u(X_{\tau_{D}}), \mathcal{J}]-\E_{x}[u(X_{\tau_{B}}), \tau_{B}=\tau_{D}]\\
&&-\E_{x}\left[\E_{X_{\tau_{B}}}[u(X_{\tau_{D}}), \mathcal{J}], \tau_{B}<\tau_{D}\right]\\
&=&v(x).
\end{eqnarray*}
\qed

\begin{thm}\label{thm:IntRep}
Let $u$ be nonnegative and harmonic on $D$ with respect to $X$. Then there exists a unique measure $\mu_{u}$ supported in $\partial D$ so that
$u(x)$ can be written as
$$
u(x)=\int_{\overline{D}^{c}}u(y)K_{D}(x,y)dy +\int_{\partial D}M_{D}(x,z)\mu_{u}(dz).
$$
\end{thm}
\pf
It follows from Theorem \hyperref[thm:Martin]{\ref{thm:Martin}} and Lemma \hyperref[lemma:rep jump3]{\ref{lemma:rep jump3}} there exists a unique measure $\mu_{u}$ supported on $\partial D$ such that
$$
u(x)-\E_{x}[u(X_{\tau_{D}}), \mathcal{J}]=\int_{\partial D}M_{D}(x,z)\mu_{u}(dz).
$$
Now it follows from \hyperref[eqn:Levy system2]{\eqref{eqn:Levy system2}} that we have
$$
u(x)=\int_{\overline{D}^{c}}u(y)K_{D}(x,y)dy +\int_{\partial D}M_{D}(x,z)\mu_{u}(dz).
$$
\qed

In \cite{Li} it is proved that there exists a bounded (classical) harmonic function on the unit disk in $\R^{2}$ that fails to have tangential limits for a.e. $\theta\in[0,2\pi]$.
Using the similar method, in \cite{Kim2, Kim} the author showed that the Stolz open sets are best possible sets for Fatou's theorem and relative Fatou's theorem for
transient censored stable processes and stable processes, respectively for $d=2$ and $D=B(0,1)$.

A curve $C_{0}$ is called a tangential curve in $B(0,1)$ if $C_{0}\cap \partial B(0,1)=\{w\}\in \partial B(0,1)$, $C_{0}\setminus \{w\}\subset B(0,1)$, and for any $r>0$ and $\beta>1$ $C_{0}\cap B(w,r)\nsubseteq A_{w}^{\beta}\cap B(w,r)$.
Let $C_{\theta}$ be a rotated curve $C_{0}$ about the origin through an angle $\theta$.
We will adapt arguments in \cite{Kim2, Kim, Li} to prove that the Stolz open sets are best possible sets for Fatou's theorem for $X$ by showing that there exists bounded harmonic function $u(x)$ with respect to $X^{B(0,1)}$ such that the tangential limit $\displaystyle\lim_{x\in C_{\theta}, x\rightarrow z}u(x)$ does not exist, where $C_{\theta}$ is a tangential curve inside $B(0,1)$.

We start with a simple lemma that is analogue to \cite[Lemma 2]{Li} (see also \cite[Lemma 3.19]{Kim2} and \cite[Lemma 3.22]{Kim}).
Let $D=B(0,1)\in \R^{2}$, $x_{0}=0$, and $\sigma_{1}$ be the normalized surface measure of $\partial B(0,1)$.
Define $h_{1}(x):=\P(X_{\tau_{B(0,1)}}\in \partial B(0,1))$ and $h_{2}(x):=\int_{\partial B(0,1)}M_{B(0,1)}(x,z)\sigma_{1}(dz)$. 
It follows from Theorem \ref{thm:Rep harmonic boundary} 
$$
H(z)=\lim_{A_{z}^{\beta}\ni x\rightarrow z\in\partial B(0,1)}\frac{h_{1}(x)}{h_{2}(x)}
$$
exists, $0<c_{1}\leq H(z)\leq c_{2}<\infty$ for some constants $c_{1},c_{2}>0$, and 
$$
h_{1}(x)=\int_{\partial B(0,1)}M_{B(0,1)}(x,z)H(z)\sigma_{1}(dz).$$
\begin{lemma}\label{lemma:near 1}
Let $h_{1}(x)=\int_{\partial B(0,1)}M_{B(0,1)}(x,z)H(z)\sigma_{1}(dz)$ 
and $U(z)$ be a nonnegative and measurable function on $\partial B(0,1)$, and $0\leq U(e^{i\theta})\leq 1$, $\theta\in [0,2\pi]$.
Suppose that $U(e^{i\theta})=1$ for $\theta_{0}-\lambda\leq\theta\leq \theta_{0}+\lambda$ for some $0<\lambda<\pi$.
Let $u(x)=\int_{\partial B(0,1)}M_{B(0,1)}(x,z)U(z)H(z)\sigma_{1}(dz)$, $x\in B(0,1)$.
Then for any $\eps>0$ there exists
$\delta=\delta(\eps,\phi)$, independent of $\lambda$, such that
$$
1-\eps \leq \frac{u(\rho e^{i\theta_{0}})}{h_{1}(\rho e^{i\theta_{0}})}\leq 1, \quad \text{if } \rho>1-\lambda\delta.
$$
\end{lemma}
\pf
Since $0\leq U(z)\leq 1$ we have
\begin{eqnarray*}
0&\leq & \frac{u(x)}{h_{1}(x)}=\frac{1}{h_{1}(x)}\int_{\partial B(0,1)}M_{B(0,1)}(x,z)U(z)H(z)\sigma_{1}(dz)\\
&\leq& \frac{1}{h_{1}(x)}\int_{\partial B(0,1)}M_{B(0,1)}(x,z)H(z)\sigma_{1}(dz)=1.
\end{eqnarray*}
Let $V(z):=\frac{1-U(z)}{2}$ so that $0\leq V(z)\leq \frac12$ and $V(e^{i\theta})=0$ for $\theta_{0}-\lambda\leq\theta\leq \theta_{0}+\lambda$.
By the triangular inequality we have $|e^{i\theta_{0}}-e^{i\theta}|\leq |e^{i\theta_{0}}-\rho e^{i\theta_{0}}|+
|\rho e^{i\theta_{0}}- e^{i\theta}|= (1-\rho)+|\rho e^{i\theta_{0}}-e^{i\theta}|$. Hence
\begin{eqnarray*}
|\rho e^{i\theta_{0}}- e^{i\theta}|&\geq& |e^{i\theta_{0}}-e^{i\theta}|-(1-\rho)\geq 2\left|\sin(\frac{\theta_{0}-\theta}{2})\right|-\delta|\theta_{0}-\theta|\\
&\geq&\frac{2}{\pi}|\theta_{0}-\theta|-\delta|\theta_{0}-\theta|\\
&=&(\frac{2}{\pi}-\delta)|\theta_{0}-\theta|
\end{eqnarray*}
for $|\theta_{0}-\theta|>\lambda$.
Hence from \hyperref[eqn:Martin]{\eqref{eqn:Martin}} we have for $\rho>1-\lambda\delta$
\begin{eqnarray*}
&&\int_{0}^{2\pi}M_{B(0,1)}(\rho e^{i\theta_{0}},e^{i\theta})V(e^{i\theta})d\theta\\
&\leq&c_{1}(1-\rho)\int_{0}^{2\pi}\frac{V(e^{i\theta})}{|\rho e^{i\theta_{0}}-e^{i\theta}|^{2}}d\theta\\
&\leq&c_{1}(1-\rho)(\frac{2}{\pi}-\delta)^{-2}\int_{|\theta-\theta_{0}|>\lambda}\frac{d\theta}{|\theta_{0}-\theta|^{2}}\\
&\leq&c_{1}\frac{1-\rho}{\lambda}(\frac{2}{\pi}-\delta)^{-2}\\
&\leq&c_{1}\frac{\delta}{(\frac{2}{\pi}-\delta)^{2}}.
\end{eqnarray*}
From Theorem \ref{thm:Rep harmonic boundary} $H(e^{i\theta})\leq c_{2}$ for some constant $c_{2}>0$.
Hence if $\delta\leq\frac{1}{\pi}$ we have
\begin{eqnarray*}
\frac{u(\rho e^{i\theta_{0}})}{h_{1}(\rho e^{i\theta_{0}})}&=&\frac{1}{h_{1}(\rho e^{i\theta_{0}})}\frac{1}{2\pi}\int_{0}^{2\pi}M_{B(0,1)}(\rho e^{i\theta_{0}},e^{i\theta})(1-2V(e^{i\theta}))H(e^{i\theta})d\theta\\
&\geq&\frac{1}{h_{1}(\rho e^{i\theta_{0}})}(h_{1}(\rho e^{i\theta_{0}})-2c_{1}c_{2}\frac{\delta}{(\frac{2}{\pi}-\delta)^{2}})\\
&\geq&1-c_{3}\delta.
\end{eqnarray*}
Now for given $\eps$ take $\delta=\frac{1}{\pi}\wedge \frac{\eps}{c_{3}}$ and we reach the conclusion of the lemma.
\qed

Once we have Lemma \hyperref[lemma:near 1]{\ref{lemma:near 1}} by adapting the argument in \cite{Li} we have the following theorem.

\begin{thm}
There exists a bounded and nonnegative harmonic function $u(x)$ with respect to
$X^{B(0,1)}$ such that for a.e. $\theta\in [0,2\pi]$ with respect to Lebesgue measure,
$$
\lim_{|x|\rightarrow 1, x\in C_{\theta}}u(x) \text{ does not exist}.
$$
\end{thm}
\pf
Let $h_{1}(x)=\P_{x}\left(X_{\tau_{B(0,1)}}\in \partial B(0,1)\right)$ as in Lemma \ref{lemma:near 1}.
By following the argument in \cite{Li} there exist nonnegative harmonic functions $u_{k}(x)$ with respect to $X^{B(0,1)}$
defined on some $E_{k}^{*}$ such that
$$
\lim_{x\rightarrow w\in \partial B(0,1)}\frac{u_{k}(x)}{h_{1}(x)}=0 \text{ radially and } \limsup_{x\rightarrow w\in \partial B(0,1)}\frac{u_{k}(x)}{h_{1}(x)}=2^{-k} \text{ along one branch of } C_{\theta}.
$$
Let $u(x)=\displaystyle\sum_{k=1}^{\infty}u_{k}(x)$. For this $u(x)$ by following the argument in \cite{Li} with Lemma \hyperref[lemma:near 1]{\ref{lemma:near 1}}
(see also \cite[Theorem 3.23]{Kim}) we have
$$
\lim_{|x|\rightarrow 1, x\in C_{\theta}}\frac{u(x)}{h_{1}(x)} \text{ does not exist for a.e. } \theta\in [0,2\pi].
$$
It follows from Proposition \hyperref[prop:limit=1]{\ref{prop:limit=1}} $\displaystyle\lim_{x\rightarrow z\in \partial B(0,1)}h_{1}(x)=1$.
Hence we have
$$
\lim_{|x|\rightarrow 1, x\in C_{\theta}}u(x)  \text{ does not exist for a.e. } \theta\in [0,2\pi].
$$
\qed

\medskip
{\bf Acknowledgement}
The author thanks Renming Song for very helpful comments and suggestions while this work was under progress and Yunju Lee whose suggestion improved many arguments in Section \hyperref[Relative Fatou theorem]{\ref{Relative Fatou theorem}}. 
The author also thanks the anonymous referee for the careful reading of the first version of this paper whose suggestions improved the quality of the paper.
\bibliographystyle{plain}
\bibliography{FTSBM}
\end{doublespace}
\vskip 0.3truein

{\bf Hyunchul Park}

Department of Mathematics, State University of New York at New Paltz, NY 12561,
USA

E-mail: \texttt{parkh@newpaltz.edu}
\end{document}